\documentclass{amsart}
\usepackage{amssymb}
\usepackage{amsfonts}
\usepackage{amsbsy}

\usepackage{latexsym}
\usepackage[mathscr]{eucal}
\usepackage{verbatim}
\usepackage[normalem]{ulem}

\theoremstyle{plain}
\newtheorem{theorem}{Theorem}[section]
\newtheorem{corollary}[theorem]{Corollary}
\newtheorem{lemma}[theorem]{Lemma}
\newtheorem{proposition}[theorem]{Proposition}
\newtheorem{remark}[theorem]{Remark}
\newtheorem{definition}[theorem]{Definition}

\newcommand{\A}{\mathcal{A}}
\newcommand{\B}{\mathcal{B}}

\newcommand{\Z}{\mathcal{Z}}
\newcommand{\K}{\mathcal{K}}
\newcommand{\eps}{\epsilon}

\DeclareMathOperator{\TA} {\mathcal{T}(\A)}

\DeclareMathOperator{\Ext} {\partial_{e}(\TA)}
\DeclareMathOperator{\her} {her}
\DeclareMathOperator{\Her} {Her}

\DeclareMathOperator{\St} {\A\otimes \K}
\DeclareMathOperator{\M}{\mathcal M(\mathcal A \otimes \mathcal K)}
\DeclareMathOperator{\Mul}{\mathcal M}
\DeclareMathOperator{\Ped}{Ped(\A)}
\DeclareMathOperator{\lin}{span}
\DeclareMathOperator{\co}{co}

\def\sideremark#1{\ifvmode\leavevmode\fi\vadjust{\vbox to0pt{\vss
\hbox to 0pt{\hskip\hsize\hskip1em
\vbox{\hsize2cm\tiny\raggedright\pretolerance10000
\noindent#1\hfill}\hss}\vbox to8pt{\vfil}\vss}}}

\newcommand{\be}{\begin{equation}\label}
\newcommand{\ee}{\end{equation}}
\newcommand{\bq}{\begin{equation*}}
\newcommand{\eq}{\end{equation*}}
\newcommand{\ba}{\begin{align*}}
\newcommand{\ea}{\end{align*}}
\newcommand{\bp}{\begin{proof}}
\newcommand{\ep}{\end{proof}}
\newcommand{\bL}{\begin{lemma}\label}
\newcommand{\eL}{\end{lemma}}
\newcommand{\bP}{\begin{proposition}\label}
\newcommand{\eP}{\end{proposition}}
\newcommand{\bC}{\begin{corollary}\label}
\newcommand{\eC}{\end{corollary}}
\newcommand{\bT}{\begin{theorem}\label}
\newcommand{\eT}{\end{theorem}}
\newcommand{\bR}{\begin{remark}\label}
\newcommand{\eR}{\end{remark}}
\newcommand{\bD}{\begin{definition}\label}
\newcommand{\eD}{\end{definition}}

\numberwithin{equation}{section}

\thanks{This work was partially supported by the Simons Foundation (grant No 245660 to Victor Kaftal and grant No  281966 to Shuang Zhang)}
\author{Victor Kaftal}

\address{Department of Mathematics\\
University of Cincinnati\\
P. O. Box 210025\\
Cincinnati, OH\\
45221-0025\\
USA}

\email{victor.kaftal@uc.edu}

\author{P. W. Ng}

\address{Department of Mathematics\\
University of Louisiana\\
217 Maxim D. Doucet Hall\\
P.O. Box 41010\\
Lafayette, Louisiana\\
70504-1010\\
USA}

\email{png@louisiana.edu}

\author{Shuang Zhang}

\address{Department of Mathematics\\
University of Cincinnati\\
P.O. Box 210025\\
Cincinnati, OH\\
45221-0025\\
USA}

\email{shuang.zhang@uc.edu}

\date{\today}


\begin{document}
\title{Strict comparison of positive elements in multiplier algebras}
\begin{abstract}

Main result: If a C*-algebra $\A$ is simple, $\sigma$-unital, has finitely many extremal traces, and has strict comparison of positive elements by traces, then its multiplier $\Mul(\A)$ also has strict comparison of positive elements by traces. The same results holds if ``finitely many extremal traces" is replaced by ``quasicontinuous scale".

A key ingredient in the proof is that every positive element in the multiplier algebra of an arbitrary $\sigma$-unital C*-algebra can be approximated  by a bi-diagonal series.
An application of strict comparison: If $\A$ is a simple separable stable $\sigma$-unital C*-algebra with real rank zero, stable rank one, strict comparison of positive elements by traces, then whether a positive element is a linear combination of projections depends on the trace values of its range projection. 

\end{abstract}
\maketitle
\
\section{Introduction}
The notion of strict comparison of positive elements in a C*-algebra plays an important role in Cuntz semigroups and has attracted an increasing interest in recent years. For instance in \cite [Corollary 4.6] {RordamZ_abs} Rordam has proven that if $\A$ is unital simple exact finite and $\Z$-stable, then $\A$ has strict comparison of positive elements by traces.  It was conjectured  in  2008 by Toms and Winter   that for separable  simple nuclear C*-algebras, $\Z$-stability is equivalent to strict  comparison of positive elements by traces (see also \cite {MatuiSato}  \cite {WinterNuclear}).

All unital, simple, exact, finite , and $\Z$-stable C*-algebras have strict comparison of positive elements by traces. This large class of C*-algebras includes irrational rotation algebra, higher dimensional simple noncommutative tori, crossed products of minimal homeomorphisms on compact metric spaces with finite covering dimension, simple unital AH-algebras with bounded dimension growth, the Jiang-Su algebra and many others. 

In a previous paper we have proven that if $\A$ is  a unital  separable simple non-elementary C*-algebra with real rank zero, and has finitely many extremal traces, and strict comparison of projections by traces, then $\M$ has strict comparison of projections by traces provided that the definition is appropriately adapted to the presence of ideals in $\M$ (\cite [Theorem 3.2]{KNZCompJOT}).

The main goal this paper is to extend this result \cite [Theorem 3.2]{KNZCompJOT} to a larger class of algebras and to strict comparison of positive elements. We will prove in Theorem \ref {T:str}  that under the assumption that a simple $\sigma$-unital C*-algebra $\A$ has only finitely many extremal traces and has strict comparison of positive elements by traces, then the property of strict comparison of positive elements by traces holds also in the multiplier algebra $\Mul(\A)$.

The condition that the extremal boundary is finite can be replaced  by the weaker condition that the algebra has quasicontinuous scale (Theorem \ref {T:comp quasicont}), but cannot be removed completely. Indeed in a subsequent paper (\cite {KNZPaper2}) we  prove that  strict comparison for the multiplier algebra can fail when the extremal boundary is infinite. 

Our original motivation for obtaining strict comparison of positive elements by traces in the multiplier algebra was to apply them to the problem of decomposing positive elements into  {\it positive combination of projections} (PCP for short), namely into sums  $\sum_1^n\lambda _j p_j$ where $p_j$ are projections in $\A$, $\lambda_j$ are positive scalars, and $n$ is a finite integer.
 
In \cite {KNZPISpan} and \cite {KNZFiniteSumsVNA} we investigated the notion of PCP in the setting of purely infinite C*-algebras and W*-algebras respectively (see also \cite {KNZChina}, \cite {KNZ Mult}, \cite {KNZTorsion}).  Focusing then on finite algebras, we proved in \cite [Theorem 6.1] {KNZComm}  that if $\A$ is a simple separable stable $\sigma$-unital C*-algebra with real rank zero, stable rank one, strict comparison of projections by traces
 and has finitely many extremal traces, then $a\in \A_+$ is a PCP if and only if  $\tau(R_a)< \infty$ for all $\tau \in \TA$, where $R_a$ denotes the range projection.  A key ingredient in the proof was Brown's interpolation theorem \cite {BrownIMP}.  

When the multiplier algebra $\M$ of an  algebra $\A$ as above has real rank zero (and thus Brown's interpolation theorem is again available) a similar result holds for $\M$: a necessary and sufficient condition for $A\in (\M)_+$ to be a PCP is that  either $\tau(R_A)< \infty$ for all those $\tau\in \TA$ for which $A$ belongs to the trace ideal $I_\tau$ or $A$ is a full element.
 (\cite [Theorem 6.4]{KNZCompJOT}).

Strict comparison of positive elements permits us to obtain in Theorem \ref {T:PCP}  precisely the same result but dropping the hypothesis that $\M$ is of real rank zero.

Our paper is organized as follows.

In the preliminary section \S 2 we review the notions of tracial simplex for non-unital C*-algebras, strict comparison of of positive elements by traces, and prove some lemmas on Cuntz subequivalence which we need for this paper. In section \S 3 we obtain some technical results  on convergence of dimension functions of cut-offs of strictly converging monotone sequences (Lemma \ref {L: monot lim}).

The main tool in the proof of our main result, but possibly also of independent interest, is Theorem \ref {T:bidiag decomp} where we prove that positive elements in the multiplier algebra of an arbitrary $\sigma$-unital C*-algebra can be approximated  by a bi-diagonal series (see Definition \ref{D:bi-diagonal}). This is an extension and improvement to arbitrary $\sigma$-unital C*-algebras of the tri-diagonal form obtained previously by Elliott \cite {Elliott}  for AF algebras and by Zhang \cite [Theorem 2.2 ]{ZhangRiesz} for real rank zero C*-algebras. 

Based on the above results, in section \S 5 we present the proof of strict comparison in the multiplier algebra (Theorem \ref {T:str}) broken into a couple of technical lemmas. 

In section \S 6 we then extend this result to the case where $\A$ has a quasicontinuous scale (see Definition \ref {D:quasicontinuous scale} and Theorem \ref {T:comp quasicont}).

In section \S 7 we apply strict comparison of positive elements in the multiplier algebra to the problem of decomposing positive elements into  positive combination of projections.   The proof of Theorem \ref {T:PCP}  is based on several steps, some of which may have independent interest.  

In Proposition \ref {P:principal} we use strict comparison of positive elements  to show that every positive element $A\in \M_+$ majorizes a scalar multiple of a  projection $P$  that generates the same ideal as $A$.

A second step is the extension and reformulation of the  ``$2\times 2$" Lemma \ref  {L: 2x2new} which played a key role in obtaining PCP decompositions in purely C*-algebras and in W*-algebras (see \cite {KNZPISpan}, \cite {KNZFiniteSumsVNA}.) This lemma also provides  bounds on the number of projections  needed for a PCP decomposition.

The third essential tool is given by Lemma  \ref {L:moving}, where we show, roughly speaking, that every $\sigma$-unital hereditary sub-algebras of $\Mul(\A)$ is *-isomorphic to a hereditary subalgebra of a unital corner of the multiplier algebra.

\section{Preliminaries}
\subsection{Pedersen ideal and approximate identities}
For a simple C*-algebra $\A$ the Pedersen ideal $\Ped$ is the minimal dense  ideal of $\A$ (\cite {PedersenMeasure}, \cite {LaursenSinclair}).  It contains all the positive elements with a local unit, that is the elements $a\in \A_+$  for which there exists $b\in \A_+$ such that $ba=a$. In fact 
$$
(\Ped)_+=\{x\in \A_+\mid x\le \sum_1^n y_j\quad\text{for some $n\in \mathbb N,$ } y_j\in \A_+ ~\text{with local unit.}\}$$

Let $\B$ be a $\sigma$-unital hereditary sub-algebra of $\A$   let $h$ be a strictly positive element of $\B$ with $\|h\|=1$, and let $e_n:=f_\frac{1}{n}(h)$ where
\be{e:feps}
f_\eps(t) =\begin{cases}
 0 &\text{for } t \in [0, \frac{\eps}{1+\eps}] \\
\frac{1+\eps}{\eps^2}t - \frac{1}{\eps} &\text{for } t \in (\frac{\eps}{1+\eps},  \eps)\\
  1 &\text{for } t\in [ \eps, 1].
\end{cases}
\ee
It is well known and routine to verify that $\{e_n\}_1^\infty$ is an approximate identity of $\B$ satisfying (\ref{e:appr id})
\be{e:appr id}
e_{n+1}e_n=e_n\quad\forall\,n                     
\ee
and $e_n\in \Ped$ for all $n$.

\subsection{Traces and dimension functions}
For a simple  C*-algebra we denote by $\TA$ the collection of the (norm) lower semicontinuous  densely defined  tracial weights on $\A_+$, henceforth, traces for short. Explicitly, a trace $\tau$ 
\begin{itemize}
\item
is an additive and homogeneous map from $\A_+$   into $[0, \infty]$ (a weight);
\item satisfies the trace condition $\tau(xx^*)=\tau(x^*x)$ for all $x\in \A$; 
\item the cone  $\{x\in \A_+\mid  \tau(x)< \infty\}$ is norm dense in $\A_+$  (thus $\tau$ is also called densely finite, or semifinite);
\item satisfies the condition $\tau(x)\le \varliminf \tau(x_n)$ for $x, x_n\in \A_+$ and $\|x_n-x\|\to 0$, or equivalently,  $\tau(x)= \lim \tau(x_n)$ for $0\le  x_n\uparrow x$ in norm.
\end{itemize}
Recall that every trace is finite on $\Ped$, and hence $\tau (e_n) < \infty$ for all $\tau\in \TA$.

Using the notations in \cite {TikuisisToms}, for every $0\ne e\in\Ped_+$) set 
\be{e:normaliz} \TA_{e\mapsto 1}:=\{\tau\in \TA\mid \tau(e)=1\}.\ee Then   $\TA_{e\mapsto 1}$ is a cone base  for $\TA$ and can be viewed as a normalization (or scale) of $\TA$. 
When equipped with the topology of pointwise convergence on Ped($\A$),  $\TA_{e\mapsto 1}$ is a Choquet simplex (\cite [Proposition 3.4]{TikuisisToms}). Set
$$\Ext  \text{ to be the collection of the extreme points of }\TA.$$
We call $\Ext $ the {\it extremal boundary} of $\TA$ and its elements {\it extreme traces.}

Given two nonzero elements $e, f\in\Ped_+$, the natural one-to-one map 
$$\TA_{e\mapsto 1}\ni \tau \mapsto \frac{1}{\tau(f)}\tau\in \TA_{f\mapsto 1}$$  is a homeomorphism which maps faces onto faces and in particular, extreme points onto extreme points.  Thus the cardinality of $\Ext$ does not depend on the normalization chosen. 

To simplify notations, we will henceforth identify $\TA$ with  $\TA_{e\mapsto 1}$.   (For more details, see \cite {TikuisisToms}, \cite {ElliottRobertSantiago},  and also \cite {KNZComm}).  

If $\A$ is unital, then $\Ped=\A$ and $\TA_{1\mapsto 1}$ coincides with the usual tracial state simplex. Thus the definition of $\TA$ that we use coincides with the standard one when $\A$ is unital, and hence, by Brown's Stabilization theorem \cite {BrownStable} also when $\A$ is stable and has a nonzero projection $p$. 

Furthermore, as remarked in  \cite[5.3]{KNZComm}, by the work of F. Combes   \cite [Proposition 4.1 and Proposition 4.4] {Combes} and Ortega, Rordam, and Thiel
\cite[Proposition 5.2]{OrtegaRordamThiel}
every $\tau \in \TA$ has a unique extension, still denoted by  $ \tau$, to a lower semicontinuous (i.e., normal) tracial weight  (trace for short)  on the enveloping von Neumann algebra $\A^{**}$. 


As usual,  the dimension function $d_\tau(\cdot)$ is defined on $\Mul(\A)_+$  as $$d_\tau(A)=:\lim_n \tau(A^{1/n})\quad \forall ~A\in \Mul(\A)_+, ~\tau\in \TA.$$ As shown  in \cite[Remark 5.3]{OrtegaRordamThiel}, $d_\tau(A)=\tau(R_A)$ where $R_A$ is the range projection of $A$. In particular 
\be{e:dim cutoff}
d_\tau\Big(\big(A-\delta)_+\big) = \tau(R_{(A-\delta)_+})=\tau (\chi_{(\delta, \|A\|]}(A)) \quad \forall ~\delta\ge 0.
\ee

We will also recall that for all $0\ne A\in \Mul(\A)_+$  both the maps \be{e:hat A LAff}\TA\ni \tau\mapsto d_\tau(A)\in [0, \infty]\ee and \be{e:dim A LAff}\TA\ni \tau\mapsto \tau(A)= \hat A(\tau) \in [0, \infty]\ee are affine, lower semicontinuous, and strictly positive. 

\subsection{Cuntz subequivalence}
Let $\A$ be a C*-algebra. If $p, q$ are projections in $\A$, $p\sim q$ (resp, $p\preceq q$) denotes the Murray - von Neumann equivalence, (resp., subequivalence) that is $p=vv^*, ~q=v^*v$ for some $v\in \A$ (resp.  $p\sim p'\le q$ for some projection $p'\in \A$). 

If $a, b\in \A_+$,  $a\preceq b$  denotes the Cuntz sub-equivalence of positive elements, that $\|a-x_nbx_n^*\|\to 0$ for some sequence $x_n\in\A$.  For ease of reference we list here the following known facts (e.g., see \cite {RordamStructure}).

\bL{L:Cuntz}
Let $\A$ be a C*-algebra, $a,b\in \A_+$, $\delta>0$. Then
\item [(i)] If $a\le b$ then $a\preceq b$.
\item[(ii)] If  $\|a-b\|< \delta$ then $(a-\delta)_+\preceq b$.
\item [(iii)] If $a\preceq b$, then there is $r\in \A$ such that  $(a-\delta )_+= rbr^*$.
\item [(iv)] If $a\preceq b$, then there is  $\delta'>0$ and  $r\in \A$ such that  $(a-\delta )_+= r(b-\delta')_+r^*$.
\item [(v)] $a+b\preceq a\oplus b$.
\item [(vi)] If $a\preceq b$ then $d_\tau(a)\le d_\tau(b)$ for all $\tau\in \TA$.
\item [(vii)] $\tau(b)\le \|b\|d_\tau (b)$ and $d_\tau((b-\delta)_+)< \frac{1}{\delta}\tau(b)$ for all $\tau\in \TA$.\eL

We will need an adaptation of \cite  [Lemma 1.1]{KuchNgPer}.
\bL{L:bounds}
Let $\A$ be a C*-algebra, $a,b\in \A_+$, and $\delta>0$. If $a\preceq (b-\delta )_+$, then for every $\epsilon >0$, $(a-\epsilon )_+=xbx^*$ for some   $x\in \A$ with $\|x\|^2\le \frac{\|a\|}{\delta}$. Furthermore, $x$ can be chosen so that $xx^*\le c_1(a-\epsilon)_+$ and $x^*x\le c_2(b-\delta)_+$ for some constants $c_1$ and $c_2$. 
\eL
\bp
As a consequence of  \cite [Proposition 2.4]{RordamStructure}, there is an $s\in\A$ for which
$$
(a-\epsilon )_+= s(b-\delta )_+s^*.$$
Then $$\|s(b-\delta )_+^{1/2}\|= \|(a-\epsilon )_+\|^{1/2}\le \|a\|^{1/2}.$$ Let  $$h_{\epsilon}(t)= \begin{cases}\frac{t}{\epsilon}&t\in [0, \epsilon]\\
1&t\in [\epsilon, \|a\|]\end{cases}\quad\text{and}\quad
g_\delta(t)= \begin{cases}\frac{1}{\delta} &t\in [0, \delta]\\
\frac{1}{t}&t\in [\delta, \|b\|]\end{cases}.$$ Then both functions are continuous and
\ba \|h_\epsilon(a)\|&=1, \quad (a-\epsilon)_+= h_\epsilon(a)(a-\epsilon)_+,\\
\|g_\delta(b)\|&=\frac{1}{\delta}, \quad (b-\delta )_+ =g_\delta(b) b(b-\delta )_+.\end{align*}
Set $$x=h_\epsilon(a)s(b-\delta )_+^{1/2}g^{1/2}_\delta(b).$$ Then 
\ba
xbx^*&=h_{\epsilon}(a)s(b-\delta )_+^{1/2}g^{1/2}_\delta (b)bg^{1/2}_\delta (b)(b-\delta )_+^{1/2}s^*h_{\epsilon}(a)\\
&=h_{\epsilon}(a)s(b-\delta )_+s^*h_{\epsilon}(a)\\
&=h_{\epsilon}(a)(a-\epsilon )_+h_{\epsilon}(a)\\
&=(a-\epsilon )_+.
 \end{align*}
Moreover,
$$
\|x\|\le\|h_{\epsilon}(a)\| \|s(b-\delta)_+^{1/2}\|\, \|g_\delta(b)^{1/2}\|\le 
\|a\|^{1/2}\, \frac{1}{\delta^{1/2}},
$$
\ba xx^*&= h_\epsilon(a)s(b-\delta )_+^{1/2}g_\delta(b)(b-\delta )_+^{1/2}s^* h_\epsilon(a)\\
&\le \frac{1}{\delta} h_\epsilon(a)s(b-\delta )_+s^* h_\epsilon(a)=\frac{1}{\delta} h_\epsilon(a)(a-\epsilon)_+ h_\epsilon(a)\\
&= \frac{1}{\delta} (a-\epsilon)_+,
\end{align*}
and 
\ba
x^*x&=g_\delta(b)^{1/2}(b-\delta)_+^{1/2}s^*h_\epsilon(a)^2s(b-\delta )_+^{1/2}g^{1/2}_\delta(b)\\
&\le \|s\|^2 g_\delta(b)^{1/2}(b-\delta)_+g^{1/2}_\delta(b)\le \frac{\|s\|^2}{\delta} (b-\delta)_+.
\end{align*}
\ep

Notice that if $a, \,b\in\A$ are selfadjoint and $a\le b$, in general it does not follow that  $a_+\le b_+$. However, we often need less. 
\bL{L:3}
Let $\A$ be a C*-algebra and $a, \,b\in\A$ be selfadjoint. 
If $a\le b$ then $a_+\preceq b_+$. 
\eL
\bp
Since  $a\le b\le b_+$ and since $\delta(t-\delta)_+\le t(t-\delta)_+$ for all $t$ and $\delta>0$, then
$$(a-\delta)_+\le \frac{(a-\delta)_+^{1/2}}{\sqrt \delta}a\frac{(a-\delta)_+^{1/2}}{\sqrt \delta}\le  \frac{(a-\delta)_+^{1/2}}{\sqrt \delta}b_+\frac{(a-\delta)_+^{1/2}}{\sqrt \delta}\preceq b_+.$$
As a consequence, $(a-\delta)_+\preceq b_+$ for all $\delta$ and hence $a_+\preceq b_+$.
\ep 
\bL{L:nifty ineq}  
Let $\A$ be a C*-algebra, $a, b\in A_+$, $\delta_i\ge 0$ with $\delta_1\ge \delta_2+\delta_3$, then
\item [(i)] $(a+b-\delta_1)_+\preceq (a-\delta_2)_+ +(b-\delta_3)_+$ 
\item [(ii)] $d_\tau(a+b)\le d_\tau(a)+d_\tau(b)$ for all $\tau\in \TA$
\item [(iii)] $d_\tau\big((a+b-\delta_1)_+\big)\le d_\tau\big((a-\delta_2)_+\big)+d_\tau\big((b-\delta_3)_+\big)$  for all $\tau\in \TA$.
\eL
\bp 
\item [(i)] Without loss of generality, $\delta_1= \delta_2+\delta_3$ Then
$$a+b-\delta_1=  (a-\delta_2)+  (b-\delta_3)\le (a-\delta_2)_++  (b-\delta_3)_+$$
hence the conclusion follows from Lemma \ref  {L:3}.
\item [(ii)] Well known
\item [(iii)] This follows from (i), the monotonicity of $d_\tau$ with respect to $\preceq$,  and (ii).
\ep

The following simple fact will be used in section \S\ref {S:2 x 2}.
\bL{L:proj prec} 
Let $\A$ be a C*-algebra, $a\in \A_+$,  $q\in \A$ be a projection, and $\delta>0$ a real number.  If $q\preceq (a-\delta)_+$, then there is a projection $p\sim q$ such that $a\ge \delta p.$  If $a\ge \delta p$ for some projection $p$ then $ p\preceq (a-\delta')_+$ for all $0\le \delta'< \delta$.  
 \eL
\bp
Assume that $q\preceq (a-\delta)_+$. Since by Lemma \ref {L:Cuntz} (v), $$\frac{1}{ 2}q= (q-\frac{1}{ 2})_+=x (a-\delta)_+x^*\quad\text{for some}\quad x\in \A,$$
it  follows that $q= (\sqrt 2 x)(a-\delta)_+(\sqrt 2 x)^*$. Thus 
$$q\sim p:= (a-\delta)_+^{\frac{1} {2}}2x^*x(a-\delta)_+^{\frac{1}{ 2}}\le 2\|x\|^2 (a-\delta)_+$$
Then $p\le R _{ (a-\delta)_+}= \chi_{(\delta, \|a\|]}(a)\le \frac{1 }{\delta}a$.
Assume now that $a\ge \delta p$ and $0\le \delta'< \delta$. then $a-\delta' \ge (\delta-\delta')p -\delta'p^\perp$ hence by Lemma \ref {L:3}
$$(\delta-\delta')p= \big((\delta-\delta')p -\delta'p^\perp\big)_+\preceq (a-\delta')_+$$
and hence $p\preceq  (a-\delta')_+.$
\ep
\subsection{Strict comparison}
When $\A$ is a simple stably finite  C*-algebra,  we say that $\A$ has  strict comparison of positive elements by traces if $\TA\ne \emptyset$ and if $a, b\in \A_+$ and $d_\tau(a)< d_\tau(b)$ for all those $\tau\in \TA$ for which $d_\tau(b)< \infty$, then  $a\preceq b$. Notice that this definition is weaker than the standard definition that requires the same property for all elements in $\St_+$ (or sometimes for all elements of $M_n(\A)_+$ and all $n\in\mathbb N$), but this weaker property is all we need  for Theorem \ref {T:str}.

Clearly, this definition would be vacuous for $\Mul(\A)$ (when $\A$ is not unital) if there is any element $b\in\A_+$ such that $d_\tau(b)=\infty$ for all $\tau$ (which is always the case when $\A$ is stable).  Indeed every element $A\in \Mul(\A)$, including full elements, would then satisfy the condition that $d_\tau(A)< d_\tau(b)$ for all those $\tau\in \TA$ for which $d_\tau(b)< \infty$. On the other hand, the condition $A\preceq b$  implies that $A$ belongs to the ideal generated by $b$, which is $\A$ by the assumption that $\A$ is simple.   To avoid this obvious obstruction we use the following definition where $I(B)$ denotes the ideal of $\Mul(\A)$ generated by $B$. 

\bD{D: strict comparison}
Let $\A$ be a simple C*-algebra with nonempty tracial simplex $\TA$. We say that $\Mul (\A)$ has strict comparison of positive elements by traces if $A\preceq B$ for $A, B\in \Mul (\A)_+$ such that $A\in I(B)$ and $d_\tau(A)< d_\tau(B)$ for all those $\tau\in \TA$ for which $d_\tau(B)< \infty$.
\eD
\subsection {Ideals in $\Mul(\A)$}
For every $\tau\in \TA$ 
$$K_\tau:=\{A\in \Mul(\A)_+\mid \tau(A)< \infty\}$$
is a hereditary cone of  $\Mul(\A)_+$ which by the trace property satisfies the condition that if  $X^*X\in K_\tau$,  then $XX^*\in K_\tau$. Let
$$L(K_\tau):=\{X\in \Mul(\A)\mid X^*X\in K_\tau\}$$
be the associated two-sided ideal of $\Mul(\A)$ and let 
$$I_\tau:= \overline {L(K_\tau)}.$$
Then it is immediate to see  from \cite [Theorem 1.5.2] {Pedersenbook} that $$I_\tau:=\overline{ \{X\in \Mul(\A)\mid \tau(X^*X)< \infty\}}=  \overline{\lin \{K_\tau\}}$$ where the closures are in norm.
The following  is also well known (for a proof see for instance \cite [Lemma 2.6] {KNZCompJOT})
\be {e: Itau} A\in (I_\tau)_+ \text{ if and only  }d_\tau\big((A-\delta)_+\big) < \infty \text{ for every  }\delta>0.\ee

\section{Dimension functions of cut-offs of monotone sequences}

\bL{L:3a}
Let $\A$ be a C*-algebra, $ T_n, T$ be normal elements of $M(\A)$, and $K \subseteq \mathbb{C}$ be a compact set for which the spectrum 
$\sigma (T_n)$ is contained in $K$ for all $n$ and 
$\sigma(T) \subseteq K$, and assume that $T_n \rightarrow T$ strictly. Then $$f(T_n) \rightarrow f(T)\quad \text{strictly}$$
for every continuous function $f : K \rightarrow \mathbb{C}$.
\eL
\begin{proof}
This is immediate when $f$ is a polynomial in one complex variable. Then apply Stone--Weierstrass Theorem.
\end{proof}
\bL {L: monot lim} Let $\A$ be $\sigma$-unital C*-algebra,  $\tau\in \TA$, $T_n, T\in \mathcal M(\A)_+$ and assume that $T_n\to T$ in the strict topology. Then
\item [(i)] If $T_n\le T_{n+1}$ for all $n$, then $d_\tau\big((T_n-\delta)_+\big) \uparrow d_\tau\big((T-\delta)_+\big)$  for all  $\delta\ge 0$.
\item [(ii)] If $T=0$, $T_n\ge T_{n+1}$ for all $n$, and $T_1\in I_\tau$, then $d_\tau\big((T_n-\delta)_+\big)\downarrow 0$ for all $ \delta >0$.
\item [(iii)] If $T_n\ge T_{n+1}$ for all $n$ and $T_1\in I_\tau$, then  for all $0< \epsilon < \delta$
$$d_\tau\big((T-\delta)_+\big)\le \lim_nd_\tau\big((T_n-\delta)_+\big) \le  d_\tau\big((T-\delta+ \epsilon)_+\big).$$

\eL

\bp Assume without loss of generality that $\|T\|\le 1$.  Since strict convergence implies strong convergence in the enveloping W*-algebra, it is easy to verify that in case (i) $T_n\le T$ and in case (ii) $T_n\ge T$. 
\item [(i)]  
Since $T_n - \delta \leq T_{n+1} - \delta \leq T - \delta$ for every  $n$ 
and hence, by Lemma \ref{L:3}, 
$$(T_n - \delta)_+ \preceq (T_{n+1} - \delta)_+ \preceq (T - \delta)_+$$
it follows that
$$d_{\tau}((T_n - \delta)_+) \leq d_{\tau}((T_{n+1} - \delta)_+) 
\leq d_{\tau}((T - \delta)_+)$$
and hence
\be{e:lim} \lim_{n} d_{\tau}((T_n - \delta)_+) 
\leq d_{\tau}((T - \delta)_+). 
\ee
Now we prove the opposite inequality. 

Since $\A$ is $\sigma$-unital, there is an approximate identity of $\A$  consisting
of an increasing sequence $e_n$ such that $e_{n+1}e_n=e_n$  for all $n$. As $T_n\to T$ strictly and since $\sigma(T), \sigma(T_n)\subset [0, 1]$ for all $n$, by Lemma \ref {L:3a}, it follows that for every $N\in \mathbb N$, 
\begin{alignat*}{2}&(T_n - \delta)_+^{1/N}\to (T - \delta)_+^{1/N} && \text{strictly}\\
&\lim _n e_k^{1/2}(T_n - \delta)_+^{1/N}e_k^{1/2}= e_k^{1/2}(T - \delta)_+^{1/N}e_k^{1/2}\hspace{2cm}&&\text{in norm}.
\end{alignat*}
Now $\tau$ is norm continuous on $e_k^{1/2} \mathcal M(\A) e_k^{1/2}= e_k^{1/2} \A e_k^{1/2} $ because $e_k\in \Ped$ and hence $\tau(e_k)< \infty$. As a consequence,   
\be {e:cont on finite}\lim _n \tau\big( e_k^{1/2}(T_n - \delta)_+^{1/N}e_k^{1/2}  \big)=\tau\big( e_k^{1/2}(T - \delta)_+^{1/N}e_k^{1/2}  \big)
\ee 
and thus
\begin{alignat*}{2}
\tau\Big( \big( (T-&\delta)_+\big)^{1/N} \Big)=\lim_k \tau\Big ( \big(T-\delta)_+\big)^{1/2N}e_k\big((T-\delta)_+\big)^{1/2N} \Big)~&(\text{normality of $\tau$})\\
&= \lim_k \tau\Big ( e_k^{1/2}\big((T-\delta)_+\big)^{1/N}e_k^{1/2}\Big)&(\text{trace property})\\
&= \lim_k \lim_n\tau\Big ( e_k^{1/2}\big((T_n-\delta)_+\big)^{1/N}e_k^{1/2}\Big)&(\text{by (\ref {e:cont on finite})}) \\
&=\lim_k \lim_n\tau\Big ( \big((T_n-\delta)_+\big)^{1/2N}e_k(T_n-\delta)_+\big)^{1/2N}\Big)&(\text{trace property}) \\
&\le \varliminf_n\tau\Big ( \big((T_n-\delta)_+\big)^{1/N}\Big)&(\text{monotonicity of $\tau$})\\
&\le \varliminf_n d_\tau\big((T_n-\delta)_+\big) &(\text{as $\|(T_n-\delta)_+\|\le 1$})\\
&=\lim_n d_\tau\big((T_n-\delta)_+\big) &(\text{as $d_\tau\big((T_n-\delta)_+\big)\uparrow$}).
\end{alignat*}
It follows that $$\lim_n d_\tau\big((T_n-\delta)_+\big)\ge \lim_N \tau\Big( \big( (T-\delta)_+\big)^{1/N} \Big)=d_\tau\big( (T-\delta)_+\big)$$ and equality follows from (\ref{e:lim}).

\item [(ii)] Let $\epsilon >0$ and let $Q_n:=\chi_{(\delta, \infty)}(T_n)$, $P_\epsilon:=\chi_{(\epsilon, \infty)}(T_1^{1/2})$. These spectral projections belong to the von Neumann algebra $\A^{**}$ and commute with $T_n$ and $T_1$ respectively.  Recall that we identify  every $\tau\in \TA$ with its extension to $\A^{**}$ (\cite[Proposition 5.2]{OrtegaRordamThiel}, see also \S 2.2)  and  that the trace of the range projection of a positive operator is just the dimension function of that operator.  In particular, 
\be{e:Qn}\tau (Q_n)= d_\tau((T_n-\delta)_+)\le d_\tau((T_1-\delta)_+)=\tau (Q_1)< \infty\ee and 
\be{e:both parts}
d_\tau((T_n-\delta)_+) \le \frac{1}{\delta}\tau(T_nQ_n) 
= \frac{1}{\delta}\Big( \tau\big(P_{\epsilon}(T_nQ_n)P_{\epsilon}\big)+\tau\big(P_{\epsilon}^\perp(T_nQ_n)P_{\epsilon}^\perp\big)\Big). 
\ee
Since  $T_1\in I_\tau$, it follows that $\tau (P_\epsilon) = d_\tau((T_1-\epsilon)_+)< \infty$ and hence $\tau$ is $\sigma$-weak continuous on $P_\epsilon\A^{**}P_\epsilon$. Therefore 
\be{e:1st part} \tau\big(P_{\epsilon}(T_nQ_n)P_{\epsilon}\big)\le  \tau\big(P_{\epsilon}T_nP_{\epsilon}\big) \to 0.\ee
Since $T_n\le T_1$, there are $G_n\in \A^{**}$ such that $T_n^{1/2}=G_nT_1^{1/2}= T_1^{1/2}G_n^*$ and $\|G_n\|\le 1$. Then 
$$
\|P_{\epsilon}^\perp T_n^{1/2}\|= \|\chi_{[0, \epsilon)}(T_1^{1/2})T_1^{1/2}G_n^*\| \le \epsilon. $$ From here and (\ref {e:Qn}) we have 
$$ \tau\big(P_{\epsilon}^\perp T_nQ_nP_{\epsilon}^\perp\big)\Big)=\tau\big(Q_nT_n^{1/2}P_{\epsilon}^\perp T_n^{1/2}Q_n\big)
\le \epsilon^2 \tau(Q_n)\le  \epsilon^2 \tau(Q_1).$$
Thus by (\ref {e:1st part}) and (\ref {e:both parts}),  it follows that $d_\tau((T_n-\delta)_+)\to 0$.
\item [(iii)] By the same argument as in part (i), $d_\tau\big ((T_n-\delta)_+\big)\downarrow$ and hence $$\lim _n d_\tau\big ((T_n-\delta)_+\big)\ge d_\tau\big ((T-\delta)_+\big).$$
By Lemma \ref {L:nifty ineq} (iii), for every $0< \epsilon < \delta$
$$d_\tau\big ((T_n-\delta)_+\big)\le d_\tau\big ((T-\delta+\epsilon)_+\big)+ d_\tau\big ((T_n-T-\epsilon)_+\big)$$
By part (ii), $\lim_n d_\tau\big ((T_n-T-\epsilon)_+\big)=0$, which concludes the proof.
\ep
\bR{R: additions}
Unlike in (i), for part  (ii) we need to assume that $\delta>0$.  Indeed, let $P$ be a projections with $0< \tau(P)< \infty$. Then $T_n:=\frac{1}{n}P\downarrow 0$ in norm, yet $d_\tau(T_n)\equiv \tau(P)\not \to  0$.\\
Similarly, in (iii) we need to assume that  $\epsilon > 0$. Indeed let as above $P$ be a projections with $0< \tau(P)< \infty$. Then $T_n:=(\delta+\frac{1}{n})P\downarrow \delta P= T$ in norm, yet $d_\tau\big ((T_n-\delta)_+\big)\equiv \tau(P)$ while $d_\tau\big ((T-\delta)_+\big)= 0$.
\eR

For ease of use in the following section, let us single out the following special case. 
\bC{C:series} Let  $\A$ be $\sigma$-unital C*-algebra, $\tau\in \TA$, and assume that
$D:=\sum_1^\infty d_k\in \mathcal M(\A)$ is the sum of a series of elements $d_k\in \A_+$ converging in the strict topology. Then
\item [(i)] $\lim_n d_\tau\Big(\big(\sum_{i=m}^n d_i-\delta\big)_+\Big)= d_\tau\Big(\big(\sum_{i=m}^\infty a_i-\delta\big)_+\Big) $ for every $\delta\ge 0$ and  $m\in \mathbb N$.
\item [(ii)] If $D\in I_\tau$  then  
$\lim_n d_{\tau}\Big(\big(\sum_{i=n}^\infty d_i-\delta\big)_+\Big)=0$ for every  $\delta >0$.
\eC

\section{Bi-diagonal decomposition}

The following theorem uses a modification of the proof for
Theorem 2.2 of \cite{ZhangRiesz}  to decompose any positive element in a general $\sigma$-unital C*-algebras into the sum of a ``bi-diagonal" series and ``small" remainder in $\A$.  Notice that there is no need to  assume the existence of an approximate identity of projections. By bi-diagonal we mean the following:
\bD{D:bi-diagonal}
Let $\A$ be a $\sigma$-unital C*-algebra and let $d_k\in\A_+$. We say that the series $D:=\sum_1^\infty d_k$ is bi-diagonal  if  
$\sum_1^\infty d_k$ converges in the strict topology
and $d_nd_m=0$ for $|n-m|\ge 2$. \eD
Every bi-diagonal series $\sum_1^\infty d_k$ can thus  written as a sum of two diagonal series, namely $\sum_1^\infty d_{2k}$, and $\sum_1^\infty d_{2k+1}$, but  the sum of two diagonal series is not necessarily bi-diagonal.
\bT{T:bidiag decomp} Let $\A$ be a $\sigma$-unital $C^*$-algebra and let  $T \in \Mul(\A)_+ $. Then for every $\epsilon
> 0$ there exist a bi-diagonal series  $\sum_1^\infty d_k$ and a self-adjoint element $a_\epsilon \in \A $
with $\| a_\epsilon\| < \epsilon$ such that 
$T = \sum_1^\infty d_k + a_\epsilon.$  Furthermore, the elements $d_k\in\Ped$ and for a fixed   approximate identity  $\{e_n\}$ of $\A$ with $e_{n+1}e_{n} = e_n$, for every $n\in \mathbb N$ there is an $N\in \mathbb N$ for which $e_n\sum_N^\infty d_k=0$.
\eT
\begin{proof}
Let $\{ e_n\}$ be an increasing approximate identity of $\A $ and as usual we assume that $e_{n+1}e_{n} = e_n $ and set $e_0=0$ (see (\ref {e:appr id})). As a consequence 
$$
(e_n-e_{n-1})(e_m-e_{m-1})=0 \quad\forall ~|n-m|\ge 2.
$$
Assume without loss of generality that $\|T\|=1$ and let $a_k:= T^{1/2}(e_k-e_{k-1})T^{1/2}$. Then $a_k\in \A_+$ for all $k$ and  $T=\sum_1^\infty a_k$ where the series converges strictly. We will construct inductively two strictly increasing sequences of positive integers $\{m_k\}_0^\infty$ and $\{n_k\}_1^\infty$ as follows.
Start by setting  $m_0=n_0=0$ and  $m_1=1$ and choosing $n_1\ge 1$ such that 
$$ \|a_1-e_{n_1}a_1e_{n_1}\|<  \frac{\epsilon}{2}\qquad\text{since $e_n\to1$ strictly and $a_1\in \A$.}$$
Now choose $m_2>m_1$ and $n_2>n_1$ such that
\begin{alignat*}{2} \|e_{n_1}\sum_{j=m_2+1}^\infty a_j \|&<\Big(\frac{\epsilon}{2^5}\Big)^2\qquad &&(\text{since $\sum_{j=m}^\infty a_j\to 0$ strictly and $e_{n_1}\in \A$})\\
\|(1-e_{n_2})\sum_{j=m_1+1}^{m_2}a_j\|&< \frac{\epsilon}{2^4}&&(\text{since $e_n\to1$ strictly and $\sum_{j=m_1+1}^{m_2}a_j\in \A$}).
\end{alignat*}
Set $b_1:=a_1$ and iterate the construction:
\begin{align} 
&\text{choose $m_k$ so that}\quad \|e_{n_{k-1}}\sum_{j=m_k+1}^\infty a_j \|< \Big(\frac{\epsilon}{2^{k+3}}\Big)^{2}\label{e:e}\\
&\text{set} \quad b_k:=\sum_{j=m_{k-1}+1}^{m_k}a_j\label{e:bk}\\
& \text{choose $n_k$ so that}\quad
\|(1-e_{n_k})b_k\|< \frac{\epsilon}{2^{k+2}}.\label{e:1-e}
\end{align}
Since $e_{n_k}$ is also an approximate unit, to simplify notations assume henceforth that $n_k=k$.
Set for all $k\ge 1$
\ba
c_1&:=e_{1}b_1e_{1}\\
c_k&:=(e_{k}-e_{k-2})b_k(e_{k}-e_{k-2})\qquad \forall \, k\ge 2.
\end{align*}
From (\ref {e:e}) (applied to $k-1$) we see that 
\ba
\|e_{k-2}b_k\|&\le \|e_{k-2}b_k^{1/2}\|=\|e_{k-2}b_ke_{k-2}\|^{1/2}\\
&\le \|e_{k-2}\sum_{j=m_{k-1}+1}^\infty
a_je_{k-2}\|^{1/2}\\ &\le \|e_{k-2}\sum_{j=m_{k-1}+1}^\infty
a_j\|^{1/2} < \frac{\epsilon}{2^{k+2}}.
\end{align*}
From the decomposition
$$
b_k-c_k= (1-e_k)b_k+ e_kb_k(1-e_k) + e_kb_ke_{k-2}+ e_{k-2}b_k(e_k-e_{k-2})
$$
and from the above inequality and (\ref {e:1-e}) we thus obtain
 that $$\|b_k-c_k\|< \frac{\epsilon}{2^k}\qquad \forall\,k.$$
As a consequence the series $a_\epsilon:=\sum_{k=1}^\infty (b_k-c_k)$ converges in norm and hence $a_\epsilon=a_\epsilon^*\in \A$. Since $ T= \sum_{k=1}^\infty a_k= \sum_{k=1}^\infty b_k$, the series $\sum_{k=1}^\infty b_k$ converges strictly and hence so does the series $D:= \sum_{k=1}^\infty c_k= T-a_\epsilon$. 
Now set $$d_k:=c_{2k-1}+c_{2k}\quad \forall \, k\ge 1.$$
Then $D= \sum_{k=1}^\infty d_k$ and 
\begin{align*}d_1&=c_{1}+c_{2} = e_1b_1e_1+e_2b_2e_2\in e_2\A e_2\\
d_k&=(e_{2k-1}-e_{2k-3})b_{2k-1}(e_{2k-1}-e_{2k-3})+ 
(e_{2k}-e_{2k-2})b_{2k}(e_{2k}-e_{2k-2}).
\end{align*} 
As a consequence,  $d_nd_m=0$ for all $|n-m|\ge 2.$

By construction, all the elements $d_k$ have a local unit and it is immediate to verify that  $ e_n\sum_N^\infty d_k=0$  for every $N\ge \frac{n+4}{2}.$
\ep

The method of the proof of Theorem \ref {T:bidiag decomp} can be applied to give a joint bi-diagonal form to multiple elements in $\Mul(\A)$  in the following sense:
\bR{R:addendum}
Let $\A$ be a $\sigma$-unital $C^*$-algebra and let  $T_1, T_2, \cdots T_N \in \Mul(\A)_+ $. Then for every $\epsilon
> 0$ there exist $N$ bi-diagonal series  $\sum_{k=1}^\infty d_{k, j}$ and self-adjoint elements $a_{\epsilon, j} \in \A $, $1\le j\le N$
with $\|a_{\epsilon, j} \| < \epsilon$ such that 
$T_j = \sum_1^\infty d_{k, j} + a_{\epsilon, j}$ 
and $d_{n, i}d_{m, j}= 0$ for $|n-m|\ge 2$ and all $1\le i,j\le N$.
\eR

Thus, up to a small remainder, every element in $\Mul(\A)_+$ is bi-diagonal and hence the sum of two diagonal series.  Diagonal series  are used extensively in multiplier algebras. We will need  the following result relating Cuntz subequivalence of (cut-offs of) summands in two diagonal series to Cuntz subequivalence of (cut-offs of) their sums. Notice that we do not need to require that the summands belong to $\A$.

\bP{P:diag Cuntz} 
Let $\A$ be a C*-algebra, $A=\sum_1^\infty A_n$, $B=\sum_1^\infty B_n$ where $A_n, B_n\in \Mul(\A)_+$, $A_nA_m=0$, $B_nB_m=0$ for $n\ne m$ and the two series converge in the strict topology. If $A_n\preceq (B_n-\delta)_+$ for some $\delta >0$ and for all $n$, then $A\preceq (B-\delta')_+$ for all $0< \delta'< \delta$. 
\eP
\bp
By Lemma \ref {L:bounds} applied to $A_n\preceq \Big((B_n-\delta')-(\delta-\delta')\Big)_+$, for every $n$ there is an $X_n\in \A$ such that 
\ba(A_n-\epsilon)_+&= X_n(B_n-\delta')X_n^*,\\ 
\|X_n\|^2&\le\frac{\|A_n\|}{\delta-\delta'}\le  \frac{\sup_n\|A_n\|}{\delta-\delta'},\\
X_nX_n^*&\le c_{1,n}(A_n-\epsilon)_+,\\
X_n^*X_n&\le c_{2,n}(B_n-\delta)_+.
\end{align*}
for some constants $c_{1,n}$ and $c_{2,n}.$ 
Therefore
\ba
 X_nX_n^*&\le \|X_n\|^2R_{(A_n-\epsilon)_+} \le \frac{\sup_n\|X_n\|^2}{\epsilon}A_n\\
X_n^*X_n&\le \|X_n\|^2R_{(B_n-\delta)_+} \le \frac{\sup_n\|X_n\|^2}{\delta}B_n.
\end{align*}
As a consequence, $R_{X_n}\le R_{A_n}$ and $R_{X_n^*}\le R_{B_n}$. But then, $$R_{X_n}R_{X_m}=R_{X_n^*}R_{X_m^*}=0 ~\forall n\ne m\quad\text{and hence}\quad 
X_nX_m^*=X_n^*X_m=0~\forall n\ne m.$$
Thus for every $m<n \in \mathbb N$ and $a\in \A$
\ba \|a\sum_{k=m}^n X_k\|^2&= \|a\sum_{k, k'=m}^n X_kX_{k'}^*a^*\|= \|a\sum_{k=m}^n X_kX_k^*a^*\|\\&\le  \frac{\sup_n\|X_n\|^2}{\epsilon}\|a\sum_{k=m}^n A_ka^*\|\le \|a\| \frac{\sup_n\|X_n\|^2}{\epsilon} \|a\sum_{k=m}^n A_k\|.
\end{align*}
Similarly
$$\|\sum_{k=m}^n X_ka\|^2\le\|a\| \frac{\sup_n\|X_n\|^2}{\delta} \|\sum_{k=m}^n B_k a\|.$$
Since  the series $\sum_1^\infty A_n$ and $\sum_1^\infty B_n$ converge strictly, it follows that $\sum_1^\infty X_n$ converges strictly. Let $X:=\sum_1^\infty X_n$.
Then $X\in \Mul(\A)$ and since $X_n=X_nR_{B_n}$  and $R_{B_n}(B_n-\delta')_+= (B_n-\delta')_+$ for every $n$, 
\ba
(A-\epsilon)_+&= \sum_1^\infty (A_n-\epsilon)_+\\
&=  \sum_1^\infty X_n(B_n-\delta')_+X_n^*\\
&= \Big(\sum_1^\infty X_n\Big)\sum_1^\infty (B_n-\delta')_+\Big(\sum_1^\infty X_n^*\Big)\\
&= X(B-\delta')_+X^*.
\end{align*}
Since $\epsilon$ is arbitrary, it follows that $A\preceq (B-\delta')_+$.
\ep
\bR{R: additions on X}
From the above proof we see that if the series $\sum_1^\infty A_n$ converges in norm, then the series $\sum_1^\infty X_n$ also converges in norm. 
\eR

\section{Strict comparison of positive elements in  $\Mul(\A)$.}
For which simple C*-algebras $\A$ does strict comparison of positive elements by traces hold for $\Mul(\A)$ when it holds for $\A$? In this section we prove that a sufficient condition is that  $\Ext$ is finite. In a subsequent  paper we will prove that when $\A$ is stable and contains a nonzero projection, finiteness of $\Ext$ is indeed necessary. In fact  when $\Ext$ is infinite even strict comparison of projections  by traces fails (\cite [Proposition 4.5]{KNZPaper2}).

We will need the following notation: for every $B\in \Mul(\A)_+$, let 
\be {e:F(B)} F(B)=\co\{ \tau\in \Ext\mid B\not \in I_{\tau}\} \ee  denote the convex combination of the extremal traces for which $B\not \in I_\tau$ and let $F(B)'$ be its complementary face (the union of the faces disjoint from $F(B)$, so the largest face disjoint from $F(B)$). Either $F(B)$ or $F(B)'$ can be empty. For this and for other basic results on convexity theory and Choquet simplexes we refer the reader to \cite {Goodearl}.
If $\Ext$ is finite, then both $F(B)$ and $F(B)'$ are closed and by \cite[Theorem 11.28]{Goodearl}, 
$$\TA=F(B)\dotplus F(B)',$$
is the direct convex sum of $F(B)$ and $F(B)'$, that is, $F(B)\cap F(B)'=\emptyset$ and every $\tau \in \TA\setminus(F(B)\cup F(B)'$ has a unique decomposition $\tau=t\mu+(1-t)\mu'$ for some $0<t<1$, $\mu\in F(B)$, and  $\mu'\in F(B)'$. Thus 
\be{e:F(B)'}F(B)'=\co\{ \tau\in \Ext\mid B\in I_{\tau}\}= \{ \tau\in \TA\mid B\in I_{\tau}\}.\ee
\bL{L: abs cont}
Let $\A$ be a  $\sigma$-unital  simple C*-algebra with finite extremal tracial boundary $\Ext$.
Let $A, B\in \Mul(\A)_+$ such that $A\in I(B)$ and $d_\tau(A)< d_\tau(B)$ for all those $\tau\in \TA$ for which $d_\tau(B)< \infty$. Then for every $\epsilon>0$ there are $\delta>0$ and $\alpha>0$  such that 
\item [(i)] $d_{\tau}((A - \epsilon)_+) +\alpha\le d_{\tau}((B - \delta)_+)<\infty$ if  $\tau\in F(B)'$
\item [(ii)] $d_{\tau}((B - \delta)_+)=\infty $ if  $\tau\in F(B).$
\eL
\bp
Notice first that if $d_\tau(B)< \infty$, then $B\in I_\tau$, i.e., $\tau\in F(B)'$.

If $d_{\tau}(B)= \infty$ for all $\tau\in \Ext$, set $\alpha:=1$, otherwise, set 
$$\alpha:=\frac{1}{2} \min \big\{ d_{\tau}(B) - d_{\tau}\Big(A\Big)\mid \tau\in \Ext, ~ d_{\tau}(B)< \infty \big\}.$$
Then it is easy to see that then
\be{e:alpha}
d_{\tau}(A)+2\alpha\le d_{\tau}(B)\quad\forall\,\tau\in \TA.
\ee
By (\ref{e: Itau}) we have
\be{e: various ineq}
d_{\tau}\big((B-\delta)_+\big)\begin{cases} < \infty &\forall \tau\in F(B)', \forall \delta >0\\
=\infty & \forall \tau\in F(B), \exists \delta >0.
\end{cases}
\ee
Since $\Ext$ is finite, we can choose $\delta_o>0$ such that $ d_{\tau}\big((B-\delta_0)_+\big)=\infty $ for all $\tau\in \Ext\cap F(B)$ and hence  for all $\tau\in F(B)$.

Since $(B-\delta)_+\uparrow B$ for $\delta\downarrow 0$ (convergence in norm), it follows that $$d_{\tau}\big((B-\delta)_+\big)\uparrow d_{\tau}(B)\quad \forall\,\tau\in \TA.$$ 
By hypothesis and (\ref {e:F(B)'}), $A\in I_{\tau}$ for all $\tau\in F(B)'$,  hence by (\ref{e: Itau}), $d_{\tau}\big((A-\epsilon)_+\big)<\infty.$
Then for  every $\tau \in F(B)'\cap \Ext$ we can choose  $0<\delta_\tau \le\delta_o$ such that  $$d_{\tau}\big((B-\delta_\tau)_+\big)\ge \begin{cases} d_\tau(B)-\alpha &\text{if }d_\tau(B)< \infty\\
d_{\tau}\big((A-\epsilon)_+\big)+\alpha &\text{if } d_\tau(B)= \infty.
\end{cases}
$$ 
Since by (\ref {e:alpha}) $$d_\tau(B)-\alpha\ge d_\tau(A)+\alpha\ge d_{\tau}\big((A-\epsilon)_+\big)+\alpha\quad \forall\, \tau\in \TA,$$
and since $\Ext$ is finite, by choosing $\delta= \min\{ \delta_\tau \mid \tau \in F(B)'\cap \Ext\}$, we have 
$$d_{\tau}\big((B-\delta)_+\big)\ge d_{\tau}\big((A-\epsilon)_+\big)+\alpha\quad \forall \,\tau \in F(B)'\cap \Ext.$$
It is then immediate to see that the same inequality holds for all $\tau \in F(B)'$. Moreover, $d_{\tau}\big((B-\delta)_+\big)< \infty$ by (\ref {e: various ineq})
which proves (i).
 Finally, (ii) follows also from (\ref {e: various ineq}) since $\delta\le \delta_o$.
 \ep

 Next is our key lemma which deals with the special case of cut-downs of bi-diagonal series. 
\bL{L:bi-diagonal case}
Let $\A$ be a  $\sigma$-unital nonunital   simple C*-algebra with strict comparison of positive elements by traces and with finite extremal tracial boundary $\Ext$, $\sum_{i=1}^\infty a_i$, $\sum_{i=1}^\infty b_i$ be two bi-diagonal series in $\Mul(\A)_+$, and let $F$ be a closed face of $\TA$ and $F'$ its complementary face  (either $F$ or $F'$  can be empty). Assume  that  
for some $\epsilon, \delta, \alpha >0$ we have  
\ba &\big(\sum\nolimits_{i=1}^\infty b_i-\delta\big)_+\not \in \A&&\\
&\sum\nolimits_1^\infty a_i\in I_{\tau}&&\text{if }\tau\in F'\\
&d_{\tau}\big((\sum\nolimits_{i=1}^\infty a_i-\epsilon)_+\big)+\alpha \le d_{\tau}\big((\sum\nolimits_{i=1}^\infty b_i-\delta)_+\big)< \infty &&\text{if }\tau\in F'\\ 
&d_{\tau}\big((\sum\nolimits_{i=1}^\infty b_i-\delta)_+\big)=\infty&&\text{if }\tau \in F.\end{align*} 
Then for any $0< \delta'< \delta$, $\epsilon'> \epsilon$
$$\Big(\sum_{i=1}^\infty a_i-2\epsilon'\Big)_+\preceq\Big (\sum_{i=1}^\infty b_i-\delta'\Big)_+.$$
\eL
\bp
The case when one of the faces $F'$ or $F$ is empty is simpler and is left to the reader, so we assume that both are non-empty.

We construct  iteratively a strictly increasing sequence of positive integers $m_k$ and two interlaced sequences of positive integers $n_k, n_k'$  with 
\be {e:interlace}n_k+2\le n_k'\le n_{k+1}-2\ee

First we notice that by Corollary  \ref {C:series} (i)$$d_{\tau}\Big(\big(\sum_{i=1}^{n} b_i-\delta\big)_+\Big)\uparrow d_{\tau}\Big(\big(\sum_{i=1}^{\infty} b_i-\delta\big)_+\Big)$$  for all $\tau\in \TA$. Since $|\Ext|< \infty$,  the convergence is uniform on $\TA$ and hence on $F'$. Thus we can choose $n_1$ so that for all $\tau\in F'$
\be{e:ineq1}
d_{\tau}\Big((\sum_{i=1}^\infty a_i-\epsilon)_+\Big)< d_{\tau}\Big(\big(\sum_{i=1}^{n_1} b_i-\delta\big)_+\Big).\ee
Since $\Big( \sum_{i=1}^\infty b_i-\delta\Big)_+\not \in \A$,  by Lemma \ref {L:nifty ineq} (i)  it follows that  for every $n\in \mathbb N$
$$\Big( \sum_{i=1}^\infty b_i-\delta\Big)_+\preceq  \sum_{i=1}^{n-1} b_i+ \Big( \sum_{i=n}^\infty b_i-\delta\Big)_+. $$ Since $\sum_{i=1}^{n-1} b_i\in \A$, it follows that
$\Big( \sum_{i=n}^\infty b_i-\delta\Big)_+\ne 0$ and hence by Lemma \ref {L:3a} $$ \Big( \sum_{i=n}^m b_i-\delta\Big)_+\to  \Big( \sum_{i=n}^\infty b_i-\delta\Big)_+$$ in the strict topology. Thus we can find  some $n_1'\ge n_1+2$ so that
\be {e:nonzero}
\Big(\sum_{i=n_1+2}^{n_1'}b_i-\delta\Big)_+\ne 0.
\ee
Now  $d_{\tau}\Big(\big(\sum_{i=m}^\infty a_i-\epsilon\big)_+\Big)\downarrow 0$ for all $\tau\in F'$ by Corollary  \ref {C:series} (ii) and again the convergence is  uniform on $F'$. Thus we can choose $m_1$ so that for all $\tau\in F'$
\be{e:ineq2}
d_{\tau}\Big(\big(\sum_{i=m_1+1}^\infty a_i-\epsilon\big)_+\Big)<d_{\tau}\Big( \big(\sum_{i=n_1+2}^{n_1'}b_i-\delta\big)_+\Big).
\ee
Finally,  $d_{\tau}\Big(\big(\sum_{i=n_1'+2}^{n}b_i-\delta\big)_+\Big)\uparrow \infty$ for all $\tau\in F$ by Corollary \ref {C:series} (i).  Furthermore, because  $|F\cap \Ext|< \infty$, we can choose $n_2\ge n_1'+2$ so that for all $\tau\in F$
\be{e:ineq3}
d_{\tau}\Big(\big(\sum_{i=1}^{m_1}a_i-\epsilon \big)_+\Big)<  d_{\tau}\Big(\big(\sum_{i=n_1'+2}^{n_2}b_i-\delta\big)_+\Big).
\ee

Set $$c_1:=\sum_{i=1}^{m_1}a_i\qquad\text{and}\qquad d_1: =\sum_{i=1}^{n_1}b_i+ 
\sum_{i=n_1'+2}^{n_2}b_i.
$$ Notice that  
\be {e:perp}\sum_{i=1}^{n_1}b_i\perp 
\sum_{i=n_1'+2}^{n_2}b_i
\ee
by the condition that $b_ib_j=0$ for $|i-j|\ge 2$.
But then for all $\tau \in F\cup F'$ we have
\begin{alignat*}{2}
d_{\tau}\big((c_1-\epsilon)_+\big)&= 
d_{\tau}\Big(\big(\sum_{i=1}^{m_1}a_i-\epsilon\big)_+\Big)\\
&< d_{\tau}\Big(\big(\sum_{i=1}^{n_1}b_i-\delta\big)_+\Big)+ d_{\tau}\Big(\big(\sum_{i=n_1'+2}^{n_2}b_i-\delta\big)_+\Big)\quad &(\text{by (\ref {e:ineq1}),  (\ref {e:ineq3})})\\
&= d_{\tau}\Big(\big(\sum_{i=1}^{n_1}b_i+\sum_{i=n_1'+2}^{n_2}b_i-\delta\big)_+\Big)&(\text{by (\ref {e:perp})})\\
&= d_{\tau}\big((d_1-\delta)_+\big).
\end{alignat*}
Since every $\tau\in \TA$ has a decomposition as a convex combination of elements in $F$ and $F'$ (unique if $\tau\not \in F\cup F'$) it follows that
\be{e:step 1 forall}
d_{\tau}\big((c_1-\epsilon)_+\big)< d_{\tau}\big((d_1-\delta)_+\big)\qquad \forall ~\tau\in \TA.
\ee

Next choose integers $n'_2, n_3$ and $m_2$ such that  $n_2+2\le n_2'\le n_3-2$ and $m_2>m_1$ and
\ba 
&\big(\sum_{i=n_2+2}^{n_2'}b_i-\delta\big)_+\ne 0&(\text{as in  (\ref{e:nonzero})})\\
&d_{\tau}\Big(\big(\sum_{i=m_2+1}^\infty a_i-\epsilon\big)+\Big)<d_{\tau}\Big( \big(\sum_{i=n_2+2}^{n_2'}b_i-\delta\big)_+\Big) \hspace{0.5 cm}  \tau\in F'\qquad &(\text{as in  (\ref{e:ineq2})})\\
&d_{\tau}\Big(\big(\sum_{i=m_1+1}^{m_2}a_i-\epsilon\big)_+\Big)<  d_{\tau}\Big(\big(\sum_{i=n_2'+2}^{n_3}b_i-\delta\big)_+\Big)\qquad  \tau\in F&(\text{as in  (\ref{e:ineq3})}).
\end{align*}
Set $$c_2:=\sum_{i=m_1+1}^{m_2}a_i\qquad\text{and}\qquad d_2: =\sum_{i=n_1+2}^{n_1'} b_i+ 
\sum_{i=n_2'+2}^{n_3}b_i.
$$ 
Because the summands come from blocks removed by more than one index we have $$\sum_{i=n_1+2}^{n_1'} b_i\perp 
\sum_{i=n_2'+2}^{n_3}b_i \quad\text{and}\quad d_1\perp d_2.$$ 
Hence as for the first step of the proof we have for all $\tau\in \TA$ 
$$d_{\tau}\big((c_2-\epsilon)_+)< d_{\tau}\Big(\big(\sum_{i=n_1+2}^{n_1'}b_i-\delta\big)_+\Big)+ d_{\tau}\Big(\big(\sum_{i=n_2'+2}^{n_3}b_i-\delta\big)_+\Big) = d_{\tau}\big((d_2-\delta)_+\big).$$
Notice that $d_2$ has a different form than the beginning term $d_1$, but from $k=3$ on we can iterate the construction keeping the same form as $d_2$. Thus we choose a strictly increasing sequence of integers $m_k$ and two interlaced sequences $n_k, n_k'$  with $$n_k+2\le n_k'\le n_{k+1}-2$$ so that setting 
$$ c_k:= \sum_{i=m_{k-1}+1}^{m_k}a_i\quad\text{and}\quad d_k:=\sum_{i=n_{k-1}+2}^{n_{k-1}'} b_i+ 
\sum_{i=n_k'+2}^{n_{k+1}}b_i.
$$
we have $\sum_{i=n_{k-1}+2}^{n_{k-1}'} b_i\perp \sum_{i=n_k'+2}^{n_{k+1}}b_i$ and
\be{e:strict always}d_{\tau}\big((c_k-\epsilon)_+\big)< d_{\tau}\big((d_k-\delta)_+\big)\qquad \forall ~k\in \mathbb N, \, \tau\in \TA.\ee
By the strict comparison of positive elements of $\A$ by traces, it follows that
\be{e: prec for summands}
(c_k-\epsilon)_+\preceq (d_k-\delta)_+\quad \forall ~k.
\ee
Moreover, from this construction we see that 
\be{e:bidiag+diag} c_ic_j=0\quad |i-j|\ge 2\quad\text{and}\quad
d_id_j=0\quad i\ne j.
\ee
Thus by construction, 
\begin{align*} 
\sum_1^\infty c_{2k}&+\sum_1^\infty c_{2k+1}=\sum_1^\infty a_i\\
\sum_1^\infty d_k&\le \sum_1^\infty b_i
\end{align*} where three series $\sum_1^\infty c_{2k}$, $\sum_1^\infty c_{2k+1}$, and $\sum_1^\infty d_k$ are diagonal (sums of mutually orthogonal terms that converge in the strict topology). 
Then by (\ref {e: prec for summands}) and Proposition \ref {P:diag Cuntz}
we obtain
\ba
\big(\sum_k c_{2k}-\epsilon'\big)_+&\preceq \Big(\sum_k d_{2k}-\delta'\Big )_+\\ \big(\sum_k c_{2k+1}-\epsilon'\big)_+&\preceq \Big(\sum_k d_{2k+1}-\delta'\Big )_+.\end{align*}

Finally, 
\begin{alignat*}{2}
(\sum_1^\infty a_i-2\epsilon')_+&\preceq \big(\sum_k c_{2k}-\epsilon'\big)_++ \big(\sum_k c_{2k+1}-\epsilon'\big)_+  &\text{(by Lemma \ref{L:nifty ineq})}\\
&\preceq \big(\sum_k c_{2k}-\epsilon'\big)_+\oplus \big(\sum_k c_{2k+1}-\epsilon'\big)_+&\text{(by Lemma \ref{L:Cuntz})}\\
&\preceq \big(\sum_k d_{2k}-\delta'\big)_++\big(\sum_k d_{2k+1}-\delta'\big)_+&\text{(since $\sum_k d_{2k}\perp \sum_k d_{2k+1}$)}\\
&=\big(\sum_k d_{k}-\delta'\big)_+&\text{(since $\sum_k d_{2k}\perp \sum_k d_{2k+1}$)}\\
&\preceq \big(\sum_i b_i-\delta'\big)_+&\text{(since $\sum_k d_{k}\le  \sum_k b_k$)}.
\end{alignat*}
\ep

Notice that in this proof we used the finiteness of $|\Ext|$ in two different ways. We used it directly to guarantee the finiteness of $F\cap \Ext$ which was essential in some steps in the proof. However the fact that  $F'\cap \Ext$ is finite was used only to guarantee that the pointwise convergence on $F'$ provided by Corollary \ref {C:series} was  uniform. Thus the conclusions of this lemma will hold also without the condition that $F'\cap \Ext$ is finite provided that the  convergence  provided by  Corollary \ref {C:series} is  uniform on $F'$. This observation will be used in the next section dealing with C*-algebras with quasi-continuous scale and so we record it as follows:

\bL{L: generalization} 
Let $\A$ be a  $\sigma$-unital nonunital  simple C*-algebra with strict comparison of positive elements by traces. Let $\sum_{i=1}^\infty a_i$ and $\sum_{i=1}^\infty b_i$ be two bi-diagonal series in $\Mul(\A)_+$. 
Let $F$ be a face $F'$ be its complementary face ($F$ or $F'$ can be empty).  Assume that 
for some $\epsilon, \delta, \alpha >0$
\item [(i)] $\big(\sum\nolimits_{i=1}^\infty b_i-\delta\big)_+\not \in \A$
\item [(ii)] $\sum\nolimits_1^\infty a_i\in I_{\tau}\quad \text{if }\tau\in F'$
\item [(iii)] $d_{\tau}\big((\sum\nolimits_{i=1}^\infty a_i-\epsilon)_+\big)+\alpha \le d_{\tau}\big((\sum\nolimits_{i=1}^\infty b_i-\delta)_+\big)< \infty \quad \text{if }\tau\in F'$
\item [(iv)] $d_{\tau}\big((\sum\nolimits_{i=1}^\infty b_i-\delta)_+\big)=\infty\quad \text{if }\tau \in F$ 
\item [(v)] $F\cap \Ext$ is finite 
\item [(vi)]  for every $m\in \mathbb N$,
$d_\tau\Big(\big(\sum_m^nb_j-\delta\big)_+\Big) \uparrow d_\tau\Big(\big(\sum_m^\infty b_j-\delta\big)_+\Big)$  uniformly on $F'$
\item [(vii)] $d_\tau\Big(\big(\sum_n^\infty a_j-\epsilon\big)_+\Big) \downarrow 0$ uniformly on $F'$. \\
Then   for any  $0< \delta'< \delta$, $\epsilon'> \epsilon$
$$\Big(\sum_{i=1}^\infty a_i-2\epsilon'\Big)_+\preceq\Big (\sum_{i=1}^\infty b_i-\delta'\Big)_+.$$
\eL

We are now in position to state and prove our main theorem.

\bT{T:str}
Let $\A$ be a  $\sigma$-unital simple C*-algebra with strict comparison of positive elements by traces and with $|\Ext|< \infty$. Then strict comparison of positive element by traces holds in $\Mul(\A)$.
\eT
\bp
Let $A, B\in \Mul(\A)_+$ such that $A\in I(B)$ and $d_\tau(A)< d_\tau(B)$ for all those $\tau\in \TA$ for which $d_\tau(B)< \infty$. Since strict comparison holds on $\A$, we can assume without loss of generality that $B\not\in \A$. Let $\epsilon >0$. By Lemma \ref {L: abs cont} we can choose $\delta >0$ and $\alpha>0$ such that 
\be{e:next reduction}\begin{cases} d_{\tau}((A - \epsilon)_+) +\alpha\le d_{\tau}((B - \delta)_+)<\infty & \tau\in F(B)'\\
d_{\tau}((B - \delta)_+)=\infty & \tau\in F(B).
\end{cases}
\ee
By Theorem \ref{T:bidiag decomp} we can find bi-diagonal decompositions  
$$A = \sum_{i=1}^\infty a_i+ a_o\quad\text{and}\quad B= \sum_{i=1}^\infty b_i+ b_o$$
where the series converge strictly and $a_i,b_i\in \A_+$, (in fact they are in $\Ped$,)   $a_ia_j=b_ib_j=0$ for $|i-j|\ge 2$, 
$a_o, b_o\in \A_{sa}$, $\|a_o\|< \epsilon$, and $\|b_o\|< \frac{\delta}{3}$. 

Since $B-b_o=\sum_{i=1}^\infty b_i\not\in \A$, and since 
$$\Big(\sum_{i=1}^\infty b_i-\delta'\Big)_+\uparrow \sum_{i=1}^\infty b_i\quad \text{for } ~\delta'\downarrow 0,$$
there is a $\delta_o>0$ such that $\Big(\sum_{i=1}^\infty b_i-\delta'\Big)_+\not \in \A$ for all $\delta'\le \delta_o$. Replacing  $\delta$ with $\min (\delta_o, \delta)$ does not decrease  $d_\tau\big((B-\delta)_+\big)$ for any $\tau\in \TA$ and assures that $(\sum_{i=1}^\infty b_i-\delta)_+\not \in \A$. 

For all  $\tau\in F(B)'$ we have $B\in I_\tau$, hence $A\in I_\tau$. Since  $a_o\in\A\subset I_\tau$,  it follows that $\sum_1^\infty a_i = A-a_o\in  I_\tau$. 

By Lemma \ref{L:Cuntz} (ii)
$$\big(\sum_{i=1}^\infty a_i-2\epsilon\big)_+\preceq (A-\epsilon)_+\quad\text{and}\quad (B-\delta)_+\preceq\big(\sum_{i=1}^\infty b_i-\frac{2\delta}{3}\big)_+.$$ Hence for all $\tau\in \TA$
$$d_\tau\Big(\big(\sum_{i=1}^\infty a_i-2\epsilon\big)_+\Big)\le d_\tau \Big((A-\epsilon)_+\Big)\quad\text{and}\quad d_\tau\Big((B-\delta)_+\Big)\le\Big(\big(\sum_{i=1}^\infty b_i-\frac{2\delta}{3}\big)_+\Big).$$

Thus 

\ba &\big(\sum\nolimits_{i=1}^\infty b_i-\delta\big)_+\not \in \A&&\\
&\sum\nolimits_1^\infty a_i\in I_{\tau}&&\text{if }\tau\in F(B)'\\
&d_{\tau}\big((\sum\nolimits_{i=1}^\infty a_i-2\epsilon)_+\big)+\alpha \le d_{\tau}\big((\sum\nolimits_{i=1}^\infty b_i-\frac{2\delta}{3})_+\big)< \infty &&\text{if }\tau\in F(B)'\\ 
&d_{\tau}\big((\sum\nolimits_{i=1}^\infty b_i-\frac{2\delta}{3})_+\big)=\infty&&\text{if }\tau \in F(B).\end{align*}
All the conditions of Lemma \ref {L:bi-diagonal case} being satisfied, we have
$$(\sum_{i=1}^\infty a_i-5\epsilon)_+\preceq (\sum_{i=1}^\infty b_i-\frac{\delta}{3})_+.$$
Since
$$(A-6\epsilon)_+\preceq(\sum_{i=1}^\infty a_i-5\epsilon)_+\quad\text{and}\quad (\sum_{i=1}^\infty b_i-\frac{\delta}{3})_+\preceq B$$
 it follows that $(A-6\epsilon)_+\preceq B$ for every $\epsilon >0$, and hence $A\preceq B$.
\ep

\section{Quasicontinuous scale}

Kucerovsky  and Perera introduced in  \cite{KuchPer} the notion of quasicontinuous scale for simple C*-algebras of real rank zero in terms of quasitraces. We adapt their definition to our setting.

\bD{D:quasicontinuous scale}
Let $\A$ be a C*-algebra  with non-empty tracial simplex $\TA$.  The lower semicontinuous affine function  $S:= \widehat {1_{\Mul(\A)}}$ is called the scale of $\A$.  \\
The scale $S$ is said to be quasicontinuous if:
\begin{itemize}
\item [(i)]  the set $F_\infty:=\{\tau\in \Ext\mid S(\tau)= \infty\}$ is finite (possibly empty);
\item [(ii)] the complementary face $F_\infty'$ of $\co( F_\infty)$ is closed (possibly empty);
\item [(iii)] the restriction $S\mid_{F_\infty'}$ of the scale $S$ to  $F_\infty'$ is continuous.
\end{itemize}
\eD
Notice that while the scale function $S$ depends on the normalization chosen  for $\TA$, the quasicontinuity of $S$ does not. Indeed,  let  $e,f$ be two positive elements in $\Ped$ and $S_e$ and $S_f$ be  the scales relative to $ \TA_{e\mapsto 1}$ and $ \TA_{f\mapsto 1}$ respectively (see \ref {e:normaliz}). Let $\psi$ be the homeomorphism $$\psi:  \TA_{e\to 1} \mapsto  \TA_{f\mapsto 1}\quad\text {given by }\psi (\tau) := \frac{1}{\tau(f)}\tau.$$
Then
$$ S_f(\psi(\tau))= \frac{S_e(\tau)}{\hat f(\tau)}\quad \forall ~\tau\in \TA_{e\mapsto 1}. $$
Since $f\in \Ped$, by the definition of the topology on $\TA$, $\hat f$ is a continuous function on $\TA$ which by the simplicity of $\A$ never vanishes, thus $ \frac{1}{\hat f(\tau)}$ is continuous. Furthermore, as stated in \S 2.2, $\psi$ maps faces onto faces, thus if $S_e$ satisfies conditions (i)-(iii) so does $S_f$.  Because of this, we can drop the reference to the specific normalization used and just refer to a scale $S$.

We can view C*-algebras with quasicontinuous scale as a natural generalization of C*-algebras with finite extremal boundary and extend to them the proof of Theorem \ref {T:str} of the previous section.

For that purpose, notice first that if $\A$ has quasicontinuous scale $S$ and  $B$ is a positive element of $ \Mul(\A)_+$, then 
the face $F(B)=\co\{ \tau\in \Ext\mid B\not \in I_{\tau}\}$ (see \ref {e:F(B)}) is contained in $\co F_\infty$ and more precisely,
\be{e:F(B) in 7}
F(B)=\co\{ \tau\in F_\infty\mid B\not \in I_{\tau}\}.
\ee
As a consequence,
\be{e:F(B)' in 7}
F(B)'= \co\{\tau\in  F_\infty\mid B\in I_\tau\}\dotplus F_\infty',
\ee  that is  $\co\{\tau\in  F_\infty\mid B\in I_\tau\}\cap F_\infty'=\emptyset$ and  every element in the face $F(B)'$ that is not in $ \co\{\tau\in  F_\infty\mid B\in I_\tau\}\cup F_\infty'$ is a unique convex combination of two elements in those faces. Both terms in this direct convex sum are closed and hence so is $F(B)'$.
 
 We start with the following lemmas. 
 \bL {L:dominated are cont}  Let $\A$ be a simple C*-algebra,  $K\subset \TA$ a closed set, and $A\le B\in \Mul(\A)_+$. If $\hat B\mid_K$ is continuous, then $\hat A\mid_K$ too is continuous.
\eL
\bp
$\hat B\mid_K= \hat A\mid_K + \widehat {B-A}\mid_K$ and since the first  function  is continuous  and the second two functions are lower semicontinuous, it is immediate to see that both must be continuous.
\ep

\bL{L:abs cont if domin} Let $\A$ be a $\sigma$-unital  simple C*-algebra, 
$K\subset \TA$ a closed set, $A, B\in \Mul(\A)_+$ with $\hat A\mid_K$ continuous, and assume that
$$d_{\tau}(A) < d_{\tau}(B)\quad\text{for all $\tau \in K$ for which $d_{\tau}(B)< \infty$.}$$
Then for every $\epsilon > 0$, there exist $\delta > 0$ and $\alpha > 0$
such that 
$$d_{\tau}((A - \epsilon)_+) +\alpha \le  d_{\tau}((B - \delta)_+)
\quad \forall ~\tau \in K.$$
If furthermore $\hat B\mid_K$ is continuous, then $d_{\tau}((B - \delta)_+)< \infty$ for all $\tau \in K.$
\eL
\bp
Assume without loss of generality that $\|A\|= 1$ and let  $f_\epsilon$ be the function defined in (\ref {e:feps}).
Then  
$$\chi_{(\epsilon, 1]}(t)\le f_{\epsilon}(t)\le \min \{1, \frac{t}{\epsilon}\},$$ 
and hence 
\begin{align}R_{(A-\epsilon)_+}&\le f_\epsilon(A),\label{e:bound 0 on f eps}\\
f_\epsilon(A)&\le \frac{1}{\epsilon}A,\label{e:bound 1 on f eps}\\  f_\epsilon(A)&\le R_A.\label{e:bound 2 on f eps}\end{align}
From (\ref {e:bound 0 on f eps}) we have
\be{e:bound 3 on f eps}
\widehat{f_\epsilon(A)}(\tau)\ge d_\tau\Big((A-\epsilon)_+\Big)\quad\forall \tau\in \TA.
\ee
From (\ref{e:bound 1 on f eps}) and Lemma \ref{L:dominated are cont} it follows that $\widehat{f_\epsilon(A)}\mid_K$ is continuous. From   (\ref{e:bound 2 on f eps}) it follows that $$\widehat{f_\epsilon(A)}(\tau)\le \tau (R_A)= d_{\tau}(A) \le  d_{\tau}(B),$$ the last inequality being strict when $d_\tau(B)< \infty$. As a consequence, the function  $\big(d_{\tau}(B)-\widehat{f_\epsilon(A)}(\tau)\big)\mid_K$ is strictly positive lower semicontinuous, and hence it has a minimum $2\alpha$ on $K$.  

Let $B_n= (B-\frac{1}{n})_+$. Then $0\le B_n\uparrow B$ (in norm) and hence $d_\tau(B_n) \uparrow d_\tau(B)$ for every $\tau\in \TA$. Since all the functions $d_\tau(B_n)$ are lower semicontinuous, by the compactness of $K$ there is an $n$ such that, 
$$d_\tau(B_n)\ge \widehat{f_\epsilon(A)}(\tau) + \alpha \quad  \forall \tau\in K.$$
Thus for $\delta:=\frac{1}{n}$ we have 
$$d_\tau((B-\delta)_+)\ge \widehat{f_\epsilon(A)}+ \alpha\ge d_\tau((A-\epsilon)_+)+\alpha \quad \forall \tau\in K,$$
where the last inequality follows from (\ref {e:bound 3 on f eps}).

If in addition $\hat B$ is continuous on $K$, then by the same reasoning as for $A$, for every $\delta>0$ we have
$$
d_\tau\big((B-\delta)_+\big)\le \widehat{f_\delta(B)}(\tau)\le \frac{1}{\delta}\hat B(\tau).
$$
Thus $d_\tau\big((B-\delta)_+\big)< \infty$ for every $\tau\in K$.
\ep

\bL {L:unif limits}  Let $\A$ be a $\sigma$-unital nonunital simple C*-algebra, $P \in \Mul(\A)$ be a projection,
$K\subset \TA$ be a closed set such that
$\widehat{P}\mid_K$ is continuous, and let\linebreak
$A:= \sum_{n=1}^{\infty} A_n$ be the strictly converging sum of elements $A_n\in (P \Mul(\A) P)_+$. Assume furthermore that there exists an approximate identity 
$\{ e_n \}_{n=1}^{\infty}$ for $P \A P$
with $e_{n+1} e_n = e_n$ for all $n \in \mathbb N$
such that for all $m \geq 1$, there exists  $N  \in \mathbb N$ with
$e_m \sum_{j = N}^{\infty} A_j = 0$.
Then for every $\delta \ge 0$, 
\item [(i)]
$d_{\tau}\Big( \Big( \sum_{j=n}^{\infty} A_j - \delta 
\Big)_+ \Big) \to 0$ uniformly on $K$.
\item [(ii)] $d_{\tau}\Big( \Big( \sum_{j=1}^n A_j - \delta \Big)_+ \Big) \to  d_{\tau}\Big( \Big( A - \delta \Big)_+ \Big)$ uniformly on $K$.
\eL
\begin{proof}

Assume without loss of generality that $\|\sum_{j = 1}^{\infty} A_j\|\le 1$ and let $\epsilon > 0$ be given.
\item [(i)]  Since  $d_{\tau}\Big( \big( \sum_{j=n}^{\infty} A_j - \delta
\big)_+ \Big) \le d_{\tau}\big( \sum_{j=n}^{\infty} A_j \big)$ for every $n$ by Lemma \ref {L:3}, it is enough to  prove the statement for $\delta=0$.

Since $e_n$ has a local unit, it belongs to the Pedersen ideal and hence by the definition of the topology on $\TA$,  $\widehat{e_n}$ is continuous. As $\widehat{e_n}\uparrow \widehat{P}$, and $ \widehat{P}\mid_K$ is continuous, by Dini's theorem the convergence is uniform on $K$. Thus choose $m$ such that $0\le \widehat{P}-\widehat{e_{m-1}}< \epsilon$  on $K$. Now choose $N$ such that 
$$e_m\sum _N^\infty A_n =0$$ 
Then for every $n\ge N$
$$\sum _n^\infty A_n= (P-e_m)\Big(\sum _n^\infty A_n\Big) (P-e_m)\le  (P-e_m)^2\le P-e_m.$$
Since $$R_{\sum _n^\infty A_n}\le R_{P-e_m}\le P-e_{m-1},$$ because $(P-e_{m-1})(P-e_m)=(P-e_m)$  we thus have for every $\tau\in K$ that
$$
 d_\tau\Big(\sum _n^\infty A_n\Big) \le\tau(P -e_{m-1})< \epsilon,$$
which proves (i).
\item [(ii)] By  Lemma \ref {L:3} and Lemma  \ref {L:nifty ineq} (iii) we have for all $n \ge 1$ and $\tau \in K$, that 
\ba
d_{\tau} \Big( \big( \sum_{j=1}^{n} A_j - \delta \big)_+ \Big)
&\leq 
d_{\tau} \Big( \big( \sum_{j=1}^{\infty} A_j - \delta  \big)_+ \Big)\\&\leq d_{\tau} \Big( \big( \sum_{j=1}^{n} A_j - \delta  \big)_+ \Big)
+ d_{\tau} \Big(\sum_{j=n+1}^{\infty} A_j\Big). 
\end{align*}
Thus (ii) follows from (i).
\end{proof}

\bR{R:tail condition}
\item [(i)] The condition that for every $n$ there exists an  $N  \in \mathbb N$ such that $e_n \sum_{j = N}^{\infty} a_j = 0$ cannot be removed. Consider for instance an element $b\in \A_+$ such that $R_b=P$ and let $a_n:=\frac{1}{2^n} f_{1/n}(b)$. Then $\sum_1^\infty a_n$ converges in norm, hence strictly, but since $R_{\sum_n^\infty a_n}=R_b$ for all $n$, it follows that $d_\tau\Big(\sum_n^\infty a_n\Big)\not\to 0.$
\item [(ii)] Substituting the continuity of $\hat P\mid _K$ with the (weaker) condition of the continuity of $\hat A\mid _K$ still permits to obtain uniform convergence on $K$ for $\delta >0$ but not for $\delta =0$. Indeed consider the case of $\A:=\B\otimes \K$ with $\B$ unital and simple, $K= \TA$, $P=1_{\Mul(\A)}$, $A_k:=  \frac{1}{2^k}1_\B\otimes e_{k,k}$, $A:=\sum_1^\infty A_k$ . Then $\hat A(\tau)=1$ for all $\tau\in \TA$, hence it is continuos, and $1_\B\otimes e_{m,m}\sum_{m+1}^\infty A_k=0$ for all $m$ but $d_\tau\big(\sum_n^\infty A_k)=\infty$ for all $\tau\in \TA$ and $n\in \mathbb N$.
\eR
\bT{T:comp quasicont}
Let $\A$ be a  $\sigma$-unital simple C*-algebra with strict comparison of positive elements by traces and with quasicontinuous scale. Then strict comparison of positive element by traces holds in $\Mul(\A)$.
\eT
\bp
In following the proof of Theorem \ref {T:str}, let  $\epsilon >0$, $A, B\in \Mul(\A)_+$ such that $B\not\in \A$, $A\in I(B)$, and $d_\tau(A)< d_\tau(B)$ for all those $\tau\in \TA$ for which $d_\tau(B)< \infty$. 

Since  $A, B\le 1_{\Mul(\A)}$ and $S\mid _{F_\infty'}$ is continuous, by 
 Lemma \ref  {L:dominated are cont}  the functions $\hat A \mid _{F_\infty'}$ and $\hat B \mid _{F_\infty'}$ are also continuous. Thus by Lemma \ref {L:abs cont if domin} there are $\delta'>0$ and $\alpha >0$ such that for all $\tau \in F'_{\infty}$
 $$d_{\tau}((A - \epsilon)_+) +\alpha \le d_{\tau}((B - \delta')_+)< \infty.$$
Moreover, for all $\tau\in F_\infty$ for which $B\in I_\tau$, we can find $\delta_\tau>0$, $\alpha_\tau>0$ such that 
$$d_{\tau}((A - \epsilon)_+) +\alpha_\tau \le d_{\tau}((B - \delta_\tau)_+)< \infty.$$
Let $\alpha= \min\{\alpha', \alpha_\tau\mid  \tau\in F_\infty, \, B\in I_\tau\}$ and $ \delta= \min\{\delta', \delta_\tau\mid  \tau\in F_\infty, \, B\in I_\tau\}.$
Since  $
F(B)'= \co\{\tau\in  F_\infty\mid B\in I_\tau\} \dotplus F_\infty'
$
(see \ref{e:F(B)'}), it thus follows that
 $$d_{\tau}((A - \epsilon)_+) +\alpha \le d_{\tau}((B - \delta)_+)< \infty \quad \forall \, \tau\in F(B)'.$$
By (\ref{e:F(B)}),
$F(B)=\co\{ \tau\in F_\infty\mid B\not \in I_{\tau}\}$
hence we can find $\delta_o>0$ such that $d_\tau\big((B-\delta_o)_+\big)=\infty$ for all $\tau\in F(B)$. By replacing if necessary $\delta$ with $\min\{\delta_o, \delta\}$, we see that the condition (\ref {e:next reduction}) in the beginning of the proof of Theorem   \ref {T:str} is satisfied. 

Thus we proceed exactly as in the proof of Theorem  \ref {T:str} and decompose $A$ and $B$ into bidiagonal series with ``small" remainders: 
$$A = \sum_{i=1}^\infty a_i+ a_o\quad\text{and}\quad B= \sum_{i=1}^\infty b_i+ b_o.$$
Recall that from Theorem \ref{T:bidiag decomp}, the bidiagonal series can be chosen so that for every $n\in \mathbb N$ there is an $N\in \mathbb N$ for which $e_n\sum_N^\infty a_k=e_n\sum_N^\infty b_k=0$ for some approximate identity $\{e_n\}$ satisfying the condition $e_{n+1}e_n=e_n$ for all $n$. Then we obtain as in the proof of Theorem  \ref {T:str}
that the conditions (i)-(iv) of Lemma \ref {L: generalization} are satisfied for $F=F(B)$ and $F'=F(B)'$.  Condition (v) holds since $F\cap \Ext\subset F_\infty$. The convergence in (vi) and (vii) is pointwise by Corollary \ref {C:series}, hence it is uniform on $ \co\{\tau\in  F_\infty\mid B\in I_\tau\}$ because this face has finite extremal boundary. In view of the decomposition in (\ref{e:F(B)'})
$F(B)'= \co\{\tau\in  F_\infty\mid B\in I_\tau\}\dotplus F_\infty'$,  we see that
 uniform convergence on $F(B)'$ holds if it holds on $F_\infty'$. Lemma \ref {L:unif limits} applied to $P=1_{\Mul(\A)}$ and $K=F_\infty'$ guarantees this uniform convergence. Thus all the condition of  Lemma \ref {L: generalization} are satisfied and the rest of the proof of Theorem  \ref {T:str} applies without changes.
 
 \ep

\section{Positive linear combinations of projections}\label{S:2 x 2}
It is well known that every element of $B(H)$ is a linear combination of projections. The same property holds for all von Neumann algebras without a  finite type  I direct summand  with infinite dimensional center \cite {GoldsteinPaszkiewicz}. However this property may fail even for C*-algebras  of real rank zero (see \cite [Proposition 5.1]{KNZComm}).

In the process of investigating linear combination of projections in C*-algebras, we found it convenient to consider the following stronger condition:

 \bD{D:univ const}
 A C*-algebra $\A$ has universal constant $V$  if every selfadjoint element $a$ in $ \A$ is a linear combination of $N$ projections  $p_j\in \A$ with $a= \sum_1^N\lambda_jp_j$ for some $N\in\mathbb N$, and $\lambda_j\in \mathbb R$, satisfying the condition
$$\sum_1^N|\lambda_j|\le V \|a\|.$$ If furthermore $N$ can be chosen independently of the element $a$, we say that $\A$ has universal constants $V$ and $N$.
 \eD

C*-algebras that have such universal constants include:
\begin{itemize}
\item von Neumann algebras without a finite type I summand with infinite dimensional center \cite{GoldsteinPaszkiewicz};
\item unital properly infinite C*-algebras (\cite [Propositions 2.6, 2.7]{KNZPISpan});
\item unital simple separable C*-algebras with real rank zero, stable rank one, strict comparison of projections and finite  extremal tracial boundary (\cite [Theorem 4.4]{KNZComm});
\item corners $P\M P$ with$P$ projection in $\M$ of  unital simple separable C*-algebras with real rank zero, stable rank one, strict comparison of projections and finite  extremal tracial boundary (\cite [Theorem 4.4] {KNZCompJOT}).
\end{itemize}

A linear combination $A=\sum_1^n\alpha_j p_j$ with projections $p_j\in \A$ and $\alpha_j>0$  will be called a {\it positive linear combination of projections} or PCP for short.

In  \cite[Proposition 2.7]{KNZPISpan} we proved that if a C*-algebra $\A$ has such universal constants and if furthermore  $\A_+$ is the closure of PCPs in $\A$ then every positive invertible element of $\A$ is a  PCP.

Thus if both conditions hold for all corners $p\A p$ of $\A$, then all positive locally invertible elements are PCP.   A key tool for the further investigation of PCP elements was the fact  that  a direct sum of projection and of  a ``small" positive perturbation is also PCP (\cite[Lemma 2.2]{KNZPISpan}). 

We can obtain the following result under less restrictive conditions:
 
\bL {L: 2x2new}
 Let $\A$ be a C*-algebra,  $p\in \A$   be a projection such that the corner algebra $p\A p$ has universal constant $V$. Then 
 \item [(i)] $p+b$ is a PCP for every $b=b^*\in p\A p$ with $\|b\|\le \frac{1}{V}$. If the corner algebra $p\A p$ has universal constants $V$ and $N$, the number of projections needed in the PCP is $N+1$.
 \item [(ii)] $p+ b$ is a PCP for every $b\in \A _+$ with $b=qb=bq$ for some projection $q\in\A$ such that $ q\perp p$ with $q\prec p$ and $\|b\|\le \frac{1}{1+V}.$ Furthermore, if $p\A p$ has 
 universal constants $V$ and $N$, $p+b$ can be decomposed as a PCP of $N+4$ projections.
  \item [(iii)]  $p+ b$ is a PCP for every $b\in \A _+$ with $b=qb=bq$ for some projection $q\in \A$ such that $  q\perp p$ with $m[q]\le [p]$ for some $m\in \mathbb N$ with $m \ge \|b\|(1+V).$ 
 \eL
\bp
 \item[(i)] By hypothesis we can find  $N$ real numbers  $\lambda_j$ and projections $q_j\in \A$ with $q_j\le p$ such that $b= \sum_{j=1}^{N}\lambda_jq_j$ and $\sum_{j=1}^{N}|\lambda_j|\le V\|b\|\le 1$.  
 $$
 p+b=  \underset{\lambda_j\ge 0}{\sum}\lambda_j q_j+  \underset{\lambda_j<0}{\sum}(-\lambda_j)(p- q_j) +\big(1+\underset{\lambda_j<0}{\sum}\lambda_j\big)p$$ 
  is a PCP of  $N+1$ projection. 
 \item [(ii)] Assume without loss of generality that $b\ne 0$ and  hence $V\|b\|< 1$ and let $\beta:=\frac{1}{1-V\|b\|}$. Then $1< \beta \le \frac{1}{\|b\|}$. Following the proof of \cite[Lemma  2.9]{KNZPISpan}, let $v\in \A$ be a partial isometry  such that $v^*v=q$ and $vv^*=p'\le p$.  Define 
 \ba
 r_1&:= \beta b + v\sqrt{\beta b- (\beta b)^2}+\sqrt{\beta b- (\beta b)^2}v^* + p'- \beta v bv^*\\
 r_2&:= \beta b - v\sqrt{\beta b- (\beta b)^2}-\sqrt{\beta b- (\beta b)^2}v ^*+ p'- \beta v bv^*
 \end{align*}
 Then $ r_1$ and $ r_2$ are projections in $\A$ and 
 $\beta b= \frac{1}{2}(r_1+r_2) -p'+ \beta v bv^*$, hence
 $$p+b=  \frac{1}{2\beta}(r_1+r_2) + \frac{1}{\beta}(p-p')+ (1-\frac{1}{\beta})\Big(p+ \frac{vbv^*}{1- \frac{1}{\beta} } \Big)$$
 Now $ 0\le \frac{vbv^*}{1- \frac{1}{\beta}} \in p'\A p'\subset p\A p$  and $\| \frac{vbv^*}{1- \frac{1}{\beta}}\|=   \frac{\|b\|}{{1- \frac{1}{\beta}}}= \frac{1}{V}$. Then by part (i) $p+ \frac{vbv^*}{1- \frac{1}{\beta} }$ is a PCP, and hence so is $p+b$. Furthermore, if $p\A p$ has 
 universal constants $V$ and $N$, by part (i),  $ p+ \frac{vbv^*}{1- \frac{1}{\beta} }$ can be decomposed as a PCP of $N+1$ projections and hence $p+b$ can be decomposed as a PCP of $N+4$ projections.
 \item [(iii)] Decompose $p=\bigoplus_{i=1}^m p_i$ into projections $p_i\in \A$ with $q\prec p_i$ for each $i$. Then 
 $$
 p+b= \sum _{i=1}^m\big(p_i+ \frac{1}{m}b\big).$$ For each $i$ it follows from part (ii) that  
 $p_i+ \frac{1}{m}b$ is a PCP and hence so is $p+b$.
 \ep
Our next  lemma permits us to embed isomorphically $\sigma$-unital hereditary sub-algebras of $\M$ into unital corners of $\M$ {\it with control on the ``size" of the corner}.  When $B\in \M$ we use the following notations: \\
$\her(B):=\overline {B(\St )B^*}$ hereditary subalgebra of $\St$\\
$\Her(B)=\overline {B\M B^*}$ hereditary subalgebra of $\M.$

\bL{L:moving}
Let $\A$ be a C*-algebra and $B\in \M_+$ be such that the hereditary algebra $\her(B)$ of $\St$, has an approximate unit $\{f_j\}$ consisting of an increasing sequence of projections. Then there is  a partial isometry $W\in (\St)^{**}$ such that
\item [(i)] $W^*W=R_B$
\item [(ii)] $WW^*\in \M$
\item [(iii)] $WB\in \M$
\item [(iv)] $W\Her(B)W^*\subseteq R\M R$ where $R:=WW^*$.
\item [(v)] The onto map 
$$\Her(B)\ni X\to \Phi(X):=WXW^*\in \Her(\Phi(B))$$
is a trace preserving  *-isomorphism of hereditary algebras.
\eL
\bp
Let $e_j:= f_j-f_{j-1}$ (with $f_0:=0$) and let $I_{\M}=\sum_1^\infty E_j$ be decomposition of the identity into projections $E_j\sim I_{\M}$. As $e_j\preceq E_j$ there are partial isometries $v_j\in \St$ such that $v_j^*v_j=e_j$ and $v_jv_j^*\le E_j$. Let $W:= \sum_1^\infty v_j.$ The series converges in the strong topology of  $(\St)^{**}$  because  both the  range projections of the partial isometries $v_j$  are mutually orthogonal and so are the range projections of $v_j^*$. Then $$W^*W=\sum_1^\infty e_j=\lim_j f_j= R_B$$ and the convergence is again in the strong topology of $(\St)^{**}$. On the other hand, 
$ WW^*= \sum_1^\infty v_jv_j^*$  in the strict topology. Thus the projection $R:=WW^*$ belongs to $ \M$. 

Next we show that $WB\in \M$. Let $a\in \St$. Then $Baa^*B\in \her(B)$, hence $f_kBa\to Ba$ in norm, or equivalently $\sum_1^n e_jBa$ converges in norm to $Ba$. Since $We_j=v_j$ for all $j$, we have
$$ W\sum_1^n e_jBa=\sum_1^n v_jBa\to WBa\in \St$$ 
since the convergence is in norm.
On the other hand, since $\sum_1^\infty v_jv_j^*$ converges strictly, $\|a\sum_n^\infty v_jv_j^*\|\to 0$ for $n\to \infty$, hence 
$$\|a\sum_n^\infty v_jv_j^*W\|= \|a\sum_n^\infty v_j\|\to0$$ and thus $aW\in \St$, whence $aWB\in \St$.

This concludes the proof of (i)-(iii).

Next, $BW^*\in \M$, hence $WB\M B W^*\subset \M$ and hence $W\Her(B)W^*\subseteq R \M R$, i.e., (iv) holds.
Finally, proving (v) is routine.
\ep

\bR{R: PZ}
\item [(i)] The above result can be seen as the construction of a projection $P\in\M$ that is equivalent to the open projection $R_B$ in the sense of Peligrad and Zsido \cite {PeligradZsido} (see also \cite {OrtegaRordamThiel}).
\item [(ii)] $\A$ has real rank zero if and only if every hereditary subalgebra of $\A$  has  an approximate identity of projections (\cite {BrownPed}). 
\eR

\bP {P:2x2 sprint}
Let $\A$ be a simple separable C*-algebra with real rank zero, 
stable rank one, strict comparison of projections, and finite extremal boundary.  Let  $P\in \M\setminus \St$ be a projection.  Then  $P+B$ is a PCP for every   $B\in (P^\perp\M P^\perp)_+$ such that  $\tau(R_B)< \infty$ for all those $\tau\in \TA$ for which $\tau(P)< \infty$. 
\eP

\bp Let $\Ext = \{\tau_j\}_1^n$ and notice that $F(B)'=\{\tau\in \TA \mid  \tau(P)< \infty\}.$  Since $\A$ has real rank zero and
$R_B$ is an open projection, it has a decomposition   $R_B=\bigoplus_1^\infty r_j$ into a strictly converging sum of mutually orthogonal projections $r_j\in \St$. By \cite [Theorem 5.1] {KNZCompJOT}, $P\M P$ has a universal constants $V$ and let $m> \|B\|(1+ V)$ be an integer.  Since $\tau( \bigoplus_1^\infty r_j)< \infty$ for all those $\tau\in \Ext$ for which $\tau(P)< \infty$ and there are only finitely many extremal traces, there exists a $k$ such that 
$\tau(\bigoplus_k^\infty r_j) < \frac{1}{m} \tau(P)$ for all $\tau\in F(B)'$. 
Let 
$$B':=B^{1/2}( \bigoplus_1^{k-1} r_j)B^{1/2}\quad\text{and} \quad B'':=B-B'= B^{1/2}( \bigoplus_k^\infty r_j)B^{1/2} $$ Then $B'\in \St_+$ and $B''\in \M_+$.  Moreover, 
$$
R_{B'}\preceq  \bigoplus_1 ^{k-1} r_j\quad{\and}\quad R_{B''}\preceq  \bigoplus_k ^\infty  r_j
$$
where  the Murray-von Neumann subequivalence $\preceq$ is in $(\St)^{**}$. Thus
$$\tau(R_{B'})< \infty \quad \forall \tau\in \Ext  \quad \text{and} \quad    \tau(R_{B''})<  \frac{1}{m} \tau(P)\quad \forall \tau\in F(B)'.$$
By \cite [Theorem 6.1] {KNZComm}, $B'$ is a PCP. Thus it remains to prove that $P+ B''$ is also a PCP.

By Lemma \ref {L:moving} there is a partial isometry $W\in (\St)^{**}$ with $W^*W= R_{B''}$, $R:=WW^*\in \M$ and that induces an isomorphism $$ \Her(B'')\ni X\to \Phi(X):=WXW^*\in  \Her(\Phi(B''))\subset R\M R.$$ Notice that $R\sim  R_{B''}$  in  $(\St)^{**}$ and hence 

$$\tau(R)=\tau(R_{B''})\le \tau (R_B)= \tau( \bigoplus_1^\infty r_j)\le \tau (P^\perp)~\,~\forall \,\tau\in \TA. $$ Thus if  $\tau(P^\perp)< \infty$, then $\tau( \bigoplus_1^\infty r_j)< \infty$, and hence $$\tau(R_{B''})\le  \tau( \bigoplus_k^\infty r_j)< \tau(P^\perp).$$ By strict comparison of projections in $\M$ (see \cite [Theorem 3.2]{KNZCompJOT}, or a consequence of Theorem \ref {T:str}) it follows that $R\preceq P^\perp$. Without loss of generality we can assume that $R\le P^\perp$. 
Now let $W':=P\oplus W$. 
Then the map
$$\Her(P\oplus B')\ni X \to \Phi'(X) = WXW^*\in \Her\big(\Phi'(P\oplus B') \big)$$
is a trace preserving *-isomorphism. 
Now $$\Phi'(P\oplus B'')= P\oplus \Phi(B'')\quad\text{and}\quad \Phi(B'')=R \Phi(B'')R.$$ Furthermore, $\tau(R)<  \frac{1}{m} \tau(P)$ for all $\tau\in F(B)'$ and since $P\not \in \St$, by the strict comparison of projections in $\M$, $m[R]\le [P].$ Then $\Phi(P\oplus B'')$ is a PCP by Lemma \ref {L: 2x2new} (iii). Since $\Phi'$ is an isomorphism of hereditary algebras, $ P+B''$ is also a PCP and hence so is $P+B$.
\ep

Next we need a result on principal ideals.  While the structure of two-sided norm closed ideals of $\Mul(\A)$ is difficult to analyze in general, a case where this structure is well understood is the following.

\bT{T:Rordam}\cite[Theorem 4.4, Proposition 4.1]{RordamIdeals}
Assume that $\A$ is a simple unital  non elementary C*-algebra with strict comparison of positive elements by traces and with finite extremal boundary  $\Ext$. Then
\item [(i)]  A proper  ideal  $\mathcal J$ of $\M$ is maximal if and only if $\mathcal J=I_\tau$ for some $\tau$ in $\Ext$.
\item [(ii)] Every proper ideal of $\M$ that properly contains $\St$ is the intersection of a (finite) number of maximal ideals. Therefore there are exactly $2^n-1$ proper ideals of $\M$ properly containing $\St$.
\eT
As a consequence if $T\in \M_+\setminus \St$ then \be {e:coroll}  I(T)= \cap \{I_\tau\mid \tau \in F(T)'\}=  \cap \{I_\tau\mid \tau \in F(T)'\cap \Ext\}. \ee

\bP{P:principal} Let $\A$ be a simple unital C*-algebra with real rank zero strict comparison of positive elements by traces and with finite extremal boundary and let $T\in \M_+\setminus \St$. Then there is a $\delta>0$ such that
\item [(i)]  $I(T)=I\big((T-\delta)_+ \big)$;
\item [(ii)] there is  a projection $P$ such that $I(P)=I(T)$ and $T\ge \delta P.$
 \eP
\bp The case when $\St=\K$ and hence $\M=B(H)$ follows from standard operator theory, so assume without loss of generality that $\A$ is not elementary. 
\item [(i)] 
By (\ref {e: Itau}), for  every $\tau\in \Ext$ for which $T\not\in I_\tau$ (that is,  for every $\tau$ in $F(T)\cap\Ext$,) there is a $\delta_\tau>0$ such that  $(T-\delta_\tau)_+\not\in I_\tau$. Let $$\delta:=\min \{\delta_\tau\mid F(T)\cap\Ext\}.$$Then $(T-\delta)_+\not\in I_\tau$ for all $\tau \in F(T)\cap\Ext$, hence $F(T)\subset F\big((T-\delta)_+\big)$. On the other hand, $(T-\delta)_+\le T\in I_\tau$ for all $\tau\in F(T)'$, hence
 $ F\big((T-\delta)_+\big)'\subset F(T)'$. Thus $ F\big((T-\delta)_+\big)'=F(T)'$ and by (\ref {e:coroll}) $I(T)=I\big((T-\delta)_+ \big)$.
 \item [(ii)] By (i), $I((T-\delta)_+) = I(T)$, and hence $$d_\tau( (T-\delta)_+)\begin{cases}< \infty&\tau\in F(T)'\\
=\infty &\tau\in F(T).
\end{cases}$$
By \cite [Proposition 3.3]{KNZCompJOT} there is a projection $Q\in \M\setminus \St$ such that $\tau(Q)=\frac{1}{2} d_\tau\big( (T-\delta)_+\big)$ for all $\tau\in \Ext$ and hence for all $\tau\in \TA$. Then $F(Q)'=F(T)'$ and hence by (\ref {e:coroll}), $I(Q)=I(T)$. Now  by strict comparison of positive elements in $\M$ (Theorem \ref {T:str}) it follows that  $Q\preceq (T-\delta)_+$. Thus by Lemma \ref {L:proj prec}, there is a projections $P\M P$ such that $T\ge \delta P$  and $P\sim Q$, and hence $I(P)=I(Q)=I(T)$.
 \ep

We list here a property we will need in the proof our the next theorem

\bL{L:pop continuity} Let $\B$ be a C*-algebra. Every $g\in C([0,1])$ is uniformly continuous on the positive part of the unit ball of $\B$.
\eL
\bp
Let $a, b$  be in the unit ball of $\B$ and let $\epsilon >0$. Find a polynomial $p_n$ such that $\|g-p_n\|_\infty< \frac{\eps}{3}$. Then $\|g(a)-p_n(a)\|< \frac{\eps}{3}$ and $\|g(b)-p_n(b)\|< \frac{\eps}{3}$. Moreover, $\|p_n(a)- p_n(b)\|\le c\|a-b\|$ where $p_n(t)= \sum_0^n \alpha_j t^j$ and $c= \sum_1^n j|\alpha_j|$. Indeed, since   $$\|a^n-b^n\|\le \|a^{n-1}\|\|a-b\|+ \|a^{n-2}\|\|b\|\|a-b\|+\cdots +\|b^{n-1}\|\|\|a-b\|\le n\|a-b\|$$ and hence $$\|p_n(a)- p_n(b)\|=\|\sum_1^n \alpha_j(a^j-b^j)\|\le \|\sum_1^n |\alpha_j|\|a^j-b^j\|\le  \|\sum_1^n j|\alpha_j|\|a-b\|$$ 
Set  $\delta= \frac{\epsilon}{3c}$. For every $\|a-b\|< \delta$ it follows that $\|p_n(a)-p_n(b)\|< \frac{\epsilon}{3}$. Thus  $\|g(a)-g(b)\|< \epsilon$. 
\ep

\bT{T:PCP}
Let $\A$ be a simple separable C*-algebra with real rank zero, stable rank one, strict comparison of projections, and finite extremal boundary, and let   $T\in \M_+$.Then $T$ is a PCP if and only if $\tau(R_T)< \infty$ for all $\tau\in F(T)'$, that is, for all $\tau$ for which $T\in I_\tau$.
\eT
\bp
We first prove the necessity. Assume that $T= \sum_{j=1}^n\lambda_jP_j$ for some $\lambda_j>0$ and projections $P_j\in \M$ and that $T\in I_\tau$ for some $\tau\in \TA$. Since $\lambda_jP_j\le T$, it follows that $P_j\in I_\tau$ and thus $\tau(P_j)<\infty.$ Let $R=\bigvee_1^n P_j\in (\St)^{**}$. Since 
$$ \tau(R_T)\le \tau(R)\le \sum_1^n\tau(P_j)$$
we conclude that $\tau(R_T)< \infty$.

Now we prove the sufficiency.
If $T\in \St$, the result is proven in \cite [Theorem 6.1]{KNZComm}. Thus assume that $T\not\in \St$ and further that $\|T\|\le 1$.
By Proposition \ref {P:principal}, there is a $0< \delta < 1$ and a projection $P\in \M$ for which $I(P)=I(T)$ and $T\ge \delta P$. Assume further that $\delta < \frac{6}{7}.$ Since $P\not \in \St$, by 
\cite {ZhangMatricial}, $P$ can be decomposed into the sum $P=P_1+P_2$ of two projections $P_1\sim P_2$. 
Then for $i=1,2$
  $$\tau(P_i) =\begin{cases} < \infty &\tau\in F(T)'\\
  =\infty &\tau\in F(T)\end{cases}$$ and hence $I(P_1)=I(P_2)=I(T)$. Set now
$$
T'= T-\frac{\delta}{ 2}P_1=T-\delta P+ \frac{\delta}{ 2}P_1+ \delta P_2.
$$
Since $T=  T'+  \frac{\delta}{ 2}P_1$, it is enough to prove that $T'$ is a PCP. Notice that $R_{T'}\le R_T$. 
 Let  $f_1$ and $f_2$ be the continuous functions defined by 
$$f_1(t) =\begin{cases} 
t &t\in [0,\frac{2}{3}\delta]\\
0&t\in [\frac{5}{6}\delta, 1]\\
\text{linear}&t\in[\frac{2}{3}\delta, \frac{5}{6}\delta]
\end{cases} \quad \text{and}\quad
f_2(t) =\begin{cases} 0 & t\in [0,\frac{2}{3}\delta]\\
t &t\in [\frac{5}{6}\delta, 1]\\
\text{linear}&t\in[\frac{2}{3}\delta, \frac{5}{6}\delta]
\end{cases}
$$
 Now consider the continuous functions $g_1$ and $g_2$ defined by 

$$g_1(t) =\begin{cases} 0 & t\in [0, \frac{\delta}{3}]\cup [\frac{2\delta}{3}, 1]\\
\frac{\delta}{2} &t=\frac{\delta}{2}\\
\text{linear}&\text{elsewhere}
\end{cases}  \quad \text{and}\quad
g_2(t) =\begin{cases} 0 & t\in [0, \frac{5\delta}{6}]\cup [\frac{7\delta}{6}, 1]\\
\delta &t=\delta\\
\text{linear}&\text{elsewhere}
\end{cases}
$$
Then for all $t\in [0,1]$, 
\begin{align} 
& f_1(t)+f_2(t)= t \label{e:sum}\\
&g_1(t)\le f_1(t) \quad \text{and}\quad g_2(t)\le f_2(t)\label{e:<=}\\
&g_1(t)f_2(t)=0 \quad \text{and}\quad g_2(t)f_1(t) =0\label{e:orthog}\\
& f_1(t)\ge \frac{\delta}{3} \quad \text{where } g_1(t)\ne 0\label{e:bound1}\\
& f_2(t)\ge \frac{5\delta}{6} \quad \text{where } g_2(t)\ne 0\label{e:bound2}.
\end{align}
Since the functions $g_1$ and $g_2$ are both continuous on $[0,1]$, by Lemma \ref {L:pop continuity} they are uniformly continuous on the set of positive contractions.
 Thus there is an integer $n$ such that $\|g_i(A)-g_i(B)\|\le \frac \delta 4$ whenever $0\le A\le 1, ~0\le B\le 1$ and $\|A-B\|\le \frac 1 n$. 

Reasoning as in the first part of the proof, we can subdivide the projections $P_1$ and $P_2$ into an orthogonal sum of $n$ projections
$$P_1= \sum_1^{n}P_{1,j}\qquad\text{and}\qquad P_2= \sum_1^{n}P_{2,j}$$
such that $I(P_{i,j})= I(T)$ for all $i=1,2$ and $1\le j\le n$. 
Then
$$
T'= \sum_{j=1}^{n}\Big(\frac{1}{n}(T-\delta P)+ \frac{\delta}{ 2}P_{1,j}+ \delta P_{2,j}\Big)
$$
Thus it is enough to prove that for every pair of projections $Q_1\perp Q_2$, with $Q_i\le R_T$,  and $I(Q_i)=I(T)$ for $i=1,2$  we have that the positive element
$$T'':= \frac{1}{n}(T-\delta P)+ \frac{\delta}{ 2}Q_1+ \delta Q_2$$
is a PCP. Notice that $R_{T''}\le R_T.$
Now $$g_1(  \frac{\delta}{ 2}Q_1+ \delta Q_2) =  \frac{\delta}{ 2}Q_1\quad\text{and}\quad g_2(  \frac{\delta}{ 2}Q_1+ \delta Q_2)=  \delta Q_2.$$ Since $\|\frac{1}{n}(T-\delta P)\|\le \frac 1 n$, it follows that
$$\|g_1(T'')- \frac{\delta}{ 2}Q_1\|= \|g_1(T'')- g_1(  \frac{\delta}{ 2}Q_1+ \delta Q_2)\|\le \frac \delta 4.$$
Then $\|\frac 2 \delta g_1(T'')- Q_1\|\le \frac 1 2$ and  by 
Lemma \ref{L:Cuntz}, $\frac 1 2 Q_1= (Q_1-\frac 1 2)_+\preceq \frac 2 \delta g_1(T'')$. Hence $Q_1\preceq  g_1(T'')$. As a consequence and by (\ref {e:bound1}), there is a projection $$Q_1'\le R_{g_1(T'')}\le \frac {1}{\delta/3}f_1(T'') \quad\text{with}\quad Q_1'\sim Q_1 \quad\text{and hence}\quad I(Q_1')=I(T).$$ Similarly, there is a projection 
$$Q_2'\le R_{g_2(T'')}\le \frac {1}{5\delta/6}f_2(T'') \quad\text{with}\quad Q_2'\sim Q_2 \quad\text{and hence}\quad I(Q_2')=I(T).$$ 

Notice that $T''=f_1(T'')+ f_2(T'')$ by (\ref {e:sum}). Then
$$
T'' = \Big (\Big (f_1(T'')- \frac {\delta}{ 3} Q'_1\Big)+ \frac{5\delta}{6}Q_2'\Big) + \Big(\Big (f_2(T'')- \frac{5\delta}{6}Q_2'\Big)+\frac {\delta}{ 3} Q'_1\Big)$$
is a decomposition of $T''$ into the sum of two positive elements.
From  (\ref {e:perp}) it follows that  $g_1(T'')f_2(T'')=0$ and hence  $Q'_2\perp R_{f_1(T'')}$.   Moreover,  $$\tau (R_{f_1(T'')})\le \tau (R_T)< \infty$$ for all $\tau\in F(T)'$ and hence for all $\tau$ for which $\tau(Q'_2)< \infty$. Similarly $Q'_1\perp R_{f_2(T'')}$ and $\tau (R_{f_2(T'')})< \infty$ for all  $\tau$ for which $\tau(Q'_1)< \infty$.
Thus both summands of $T''$  satisfy the conditions of Proposition \ref  {P:2x2 sprint} and hence are a PCP, which concludes the proof.

\ep

\end{document}